\documentclass{amsart}
\usepackage{amsmath,amssymb,amsthm,amscd}
\newtheorem{theo}{Theorem}
\newtheorem{defi}[theo]{Definition}
\newtheorem{lemm}[theo]{Lemma}
\newtheorem{prop}[theo]{Proposition}
\newtheorem{cor}[theo]{Corollary}

\theoremstyle{definition}
\newtheorem{rema}[theo]{Remark}

\begin{document}

\title[Actions of compact ableian groups on subfactors]{Classification of actions of compact abelian groups on subfactors with index less than $4$}
\author{Koichi Shimada}
\email{shimada@ms.u-tokyo.ac.jp}
\address{Department of Mathematical Sciences
University of Tokyo, Komaba, Tokyo, 153-8914, Japan}
\date{}
\begin{abstract}
We classify actions of discrete abelian groups on some inclusions of von Neumann algebras, up to cocycle conjugacy. 
As an application, we classify actions of compact abelian groups on the inclusions of AFD (approximately finite dimensional) factors of type $\mathrm{II}_1$ with index less than $4$, up to stable conjugacy. 
\end{abstract}

\maketitle

\section{Introduction}
In this paper, we classify actions of compact abelian groups on the inclusions of AFD (approximately finite dimensional) factors of type $\mathrm{II}_1$ with index less than $4$, up to stable conjugacy (Theorem \ref{cptabl}). 
This work is a natural continuation of classification of group actions on the single AFD factors. Classification of group actions was first studied by Connes \cite{C}. 
He has shown the uniqueness of the outer actions of the integer group $\mathbf{Z}$ on the AFD factor of type $\mathrm{II}_1$, up to cocycle conjugacy.
 His theorem is not only very beautiful but also quite non-trivial. 
Indeed, there are many ways of constructing outer actions of $\mathbf{Z}$ on the AFD factor of type $\mathrm{II}_1$, most of which are mutually non-conjugate. 
However, by his theorem, when we consider classification up to cocycle conjugacy, all of them are the same. 
Hence classification of group actions attracts many researcher's attention, and actions of discrete amenable groups on the AFD factors are completely classified, up to cocycle conjugacy (See Ocneanu \cite{O2} and Katayama--Sutherland--Takesaki \cite{KtST}).
 One of the next problems is to classify actions of continuous groups. For this direction, actions of compact abelian groups are completely classified, up to stable conjugacy (See Jones--Takesaki \cite{JT} and Kawahigashi--Takesaki \cite{KT}). 
In general, actions of continuous groups are not easy to classify. 
However, duals of compact abelian groups are discrete. 
By making use of this fact, they classify actions of compact abelian groups. 
Here, we consider a generalization of their results for another direction. 
We consider actions on subfactors.

Group actions on subfactors themselves have been studied by many hands. 
Constructing invariants for actions of discrete amenable groups and classifying them by these invariants have intensively been studied in 1990's and early 2000's (See Loi \cite{L}, Kawahigashi \cite{Kwh3}, Popa \cite{P2} and Toshihiko Masuda \cite{M}).  
Hence as a next problem, we consider a classification of actions of continuous groups.

Our main theorem is Theorem \ref{cptabl}, which is a complete classification of actions of compact abelian groups on the inclusions of AFD factors of type $\mathrm{II}_1$ with index less than $4$.
The theorem provides a complete list of obstructions for stable conjugacy. 
As we have explained, this result has two aspects. 
One is that this theorem is a subfactor counterpart of the results of Jones--Takesaki \cite{JT} and Kawahigashi--Takesaki \cite{KT}. 
The other is that this theorem is the first result about classification of actions of continuous groups on subfactors.

In order to classify actions of compact abelian groups on subfactors, we need to investigate actions of discrete abelian groups on inclusions of von Neumann algebras which may not be factors. 
It is possible to construct invariants by the same method as that in the case of actions of discrete amenable groups on subfactors (Section 4 of Masuda \cite{M}). 
We will try to prove the classification theorem by the same line as that of the case of single factors (Jones--Takesaki \cite{JT}, Kawahigashi--Takesaki \cite{KT}). 
In order to achieve this, we will make use of a similarity between factors of type III and subfactors. 
In the case of factors of type III, a modular automorphism group plays an important role. 
Instead, in the case of subfactors, we use a centrally trivial action of a finite abelian group. 
However, still some technical parts are different from that of the case of single factors. 
We need to handle these problems.

This paper is organized as follows. 
In Section 2, we will explain basic facts about subfactors and their automorphisms. 
We will also introduce some propositions which are necessary for our classification. 
In Section 3, we will construct the invariants for the classification. 
This part is essentially the same as Section 5 of Loi \cite{L} and Section 4 of Masuda \cite{M}. 
In Section 4, we will show a classification theorem for actions of discrete abelian groups on some inclusions of von Neumann algebras. 
The outline of the proof is the same as that of the case of single factors (Jones--Takesaki \cite{JT} or Kawahigashi--Takesaki \cite{KT}). 
However, some parts are technically different (For example, the proof of Lemma \ref{point}). 
We will handle these problems. 
In Section 5, by taking the dual, we will show the classification theorem of actions of compact abelian groups on the inclusions of AFD factors of type $\mathrm{II}_1$ with index less than $4$, up to stable conjugacy. 
In Section 6, we will further investigate actions of one dimensional torus on subfactors. 
In general, it is not easy to compute the invariants introduced in Theorem \ref{cptabl} of Section 5.
However, for actions of one dimensional torus, to some extent, it is possible to compute these invariants and the principal graphs of the crossed products. Section 6 is based on the argument of Section 6 of Loi \cite{L}.

\bigskip

\textbf{Acknowledgment}
The auther is thankful to Professor Yasuyuki Kawahigashi, who is his adviser, for his useful comments on this work. In particular, the author is thankful for his pointing out a mistake in Theorem 3. The author is also thankful to Professor Toshihiko Masuda for giving him some useful advice, particularly about Subsection 5.2. The author is supported by Research Fellowships of the Japanese Society for the Promotion of Science for Young Scientists No.26-6590. This work is also supported by the Program for Leading Graduate 
Schools, MEXT, Japan.
 
\section{Preliminaries}
\subsection{Strong amenability}
In this subsection, we explain the strong amenability defined by Popa \cite{P} 
(See subsections 2.1, 2.3, Definition 3.1.1 and Proposition 3.2.2 of Popa \cite{P}).  
Let $M\supset N$ be an inclusion of factors and $\{ M_k\}$ be the Jones tower of the inclusion $M\supset N$. 
We consider an inclusion $\mathcal{M}\supset \mathcal{N}$ with a normal faithful expectation $\mathcal{E}:\mathcal{M}\to \mathcal{N}$ and an inclusion $i:M\to \mathcal{M}$ satisfying the following commuting diagram.

$$
\begin{CD}
\mathcal{M}@>\mathcal{E}>> \mathcal{N} \\
@A i AA @A i AA \\
M @>E>> N
\end{CD}
$$

The inclusion $\mathcal{M}\supset \mathcal{N}$ is said to be a smooth representation of $M\supset N$ if the following two conditions are satisfied.

(1) We have $\overline{\mathrm{span}}M\mathcal{N}=\mathcal{M}$.

(2) We have $N'\cap M_k =\mathcal{N}'\cap M_k$ for any $k\geq 1$.
\begin{defi}
A subfactor $M \supset N$ is said to be amenable if for any smooth representation

$$
\begin{CD}
\mathcal{M}@>\mathcal{E}>> \mathcal{N} \\
@A i AA @A i AA \\
M @>E>> N,
\end{CD}
$$
there exists a \textup{(}possibly non-normal\textup{)} conditional expectation $\delta :\mathcal{M} \to M$ satisfying the following commuting diagram.

$$
\begin{CD}
\mathcal{M}@>\mathcal{E}>> \mathcal{N} \\
@V\delta VV @V \delta VV \\
M @>E>> N.
\end{CD}
$$
 
Furthermore, if the von Neumann algebra $\bigvee _k\left( M_k\cap M'\right)$ is a factor, then the inclusion $M \supset N$ is said to be strongly amenable. 
\end{defi}
Strong amenability is important for classification. Indeed, Popa \cite{P} has shown that the strongly amenable subfactors of type $\mathrm{II}_1$ are completely classified by their standard invariants (See Section 5 of Popa \cite{P}).
\subsection{Extremality}
(See Pimsner--Popa \cite{PP}). An inclusion $M\supset N$ of semifinite factors is said to be extremal if the conditional expectation with respect to the trace is minimal. 
If the index is less than $4$, then this condition is automatically satisfied. 
An inclusion $M\supset N$ of the form $A\otimes (P\supset Q)$ is said to be extremal if $P\supset Q$ is extremal, where $A$ is an abelian von Neumann algebra and $P\supset Q$ is an inclusion of semifinite factors.
 
\subsection{Automorphisms and group actions of inclusions}
In this subsection, we explain basic facts about automorphisms and group actions on inclusions of von Neumann algebras.
A basic reference is Section 4 of Loi \cite{L}. 
Let $M \supset N$ be an inclusion of von Neumann algebras. Set
\[ \mathrm{Aut}(M,N):=\{ \alpha \in \mathrm{Aut}(M)\mid \alpha (N)=N\} .\]
Then this is a closed subgroup of $\mathrm{Aut}(M)$ with respect to the u-topology, that is, $ \alpha _i \to \alpha $
if $ \| \phi \circ \alpha _i -\phi \circ \alpha \| \to 0$
for any $\phi \in M_*$. 
Let $G$ be a locally compact group.
A continuous group homomorphism from $G$ to $\mathrm{Aut}(M,N)$ is said to be an action of $G$ on the inclusion $M\supset N$. 
Two actions $\alpha $ and $\beta $ of $G$ are said to be mutually cocycle conjugate if there exist an $\alpha$-cocycle $\{ u_g\}$ of $N$ and an automorphism $\theta \in \mathrm{Aut}(M,N)$ with 
\[ \mathrm{Ad}u_g\circ \alpha _g =\theta \circ \beta _g\circ \theta ^{-1}\]
for any $g\in G$.
An automorphism $\alpha $ is said to be inner if there exists a unitary $u$ of $N$ with $\alpha =\mathrm{Ad}u$.
We denote the set of all inner automorphisms by $\mathrm{Int}(M,N)$, which is a normal subgroup of $\mathrm{Aut}(M,N)$. 
Let $\omega $ be an ultrafilter of $\mathbf{N}$. Set
\[ \mathcal{C}_\omega :=\{ (x_n)\in l^\infty (N)\mid \| x_n \phi -\phi x_n\| \to 0 \ (n\to \omega ) \ \mathrm{for} \ \mathrm{any } \ \phi \in M_*\} .\]
\[ \mathcal{I}_\omega :=\{ (x_n) \in l^\infty (N) \mid x_n \to 0 \ \mathrm{in} \ \mathrm{the} \ \mathrm{strong*} \ \mathrm{topology}\} .\]
Then $\mathcal{I}_\omega $ is a closed ideal of the $\mathrm{C}^*$-algebra $\mathcal{C}_\omega$. 
Let $C_\omega (M,N):=\mathcal{C}_\omega/\mathcal{I}_\omega$ be the quotient $\mathrm{C}^*$-algebra. 
Denote the quotient $\mathcal{C}_\omega \to C_\omega (M,N)$ by $\pi$.
Then $C_\omega (M,N)$ has a von Neumann algebra structure (See Section 4 of Loi \cite{L}).
An automorphism $\alpha $ of $M\supset N$ defines an automorphism of $C_\omega (M,N)$ by
\[ \alpha (\pi (x_n)):=\pi (\alpha (x_n)).\]
An automorphism $\alpha $ is said to be centrally trivial if $\alpha \in \mathrm{Aut}(C_\omega (M,N))$ is identity. 
We denote the set of all centrally trivial automorphisms by $\mathrm{Cnt}(M,N)$, which is a normal subgroup of $\mathrm{Aut}(M,N)$.

\subsection{Probability index}
Let $M\supset N$ be an inclusion of von Neumann algebras with a conditional expectation $E$ from $M$ to $N$. 
Then the probability index of $E $ is defined by the following (See Section 2 of Pimsner--Popa \cite{PP}).
\[ \mathrm{Ind}(E):=\mathrm{sup}\{ \lambda >0 \mid  E-\lambda ^{-1}\mathrm{Id}_M \ \mathrm{is} \ \mathrm{completely } \ \mathrm{positive}\}.\]
When $M\supset N$ is an inclusion of factors, the probability index coincides with the usual index defined by Jones \cite{J}, Kosaki \cite{Ko2} (See Theorem 2.2 of Pimsner--Popa \cite{PP} and Theorem 4.1 of Longo \cite{Longo}). 
Let $G$ be a locally compact group. 
Let $\alpha $ be an action of $G$ on the inclusion $M\supset N$. 
Then it is possible to consider the crossed product $(M\supset N) \rtimes _\alpha G$.
Let $\phi_0$ be a normal faithful state of $N$.  
Set $\phi :=\phi _0 \circ E$.
Then the modular automorphism group $\{ \sigma ^{\hat{\phi}}_t\} _{t\in \mathbf{R}}$ of the dual weight $\phi$ is of the following form.
\[ M\ni x \mapsto \sigma ^\phi _t(x) \in M,\]
\[ \lambda _s \mapsto \lambda _s[D\phi \circ \alpha _s :D\phi ]_t\]
for $s\in G$, $t\in \mathbf{R}$, where $[D\phi \circ \alpha _s:D\phi ]$ is the Connes cocycle between $\phi \circ \alpha _s$ and $\phi$. 
In general, we have $[D\phi \circ \alpha _s:D\phi ]_t\not \in N\rtimes _\alpha G$. 
However, if we assume that the action $\alpha$ commutes with the expectation $E$, then we have $[D\phi \circ \alpha _s:D\phi ]_t =[D\phi _0 \circ \alpha _s : D\phi _0]_t\in N\rtimes _\alpha G$.
Hence the action $\{\sigma ^{\hat{\phi}}_t\}$ preserves $N\rtimes _\alpha G$ globally. 
In the following, we assume that the action $\alpha $ commutes with the expectation $E$.
Under this assumption, by Takesaki's criterion, there exists a conditional expectation $\tilde{E}$ from $M\rtimes _\alpha G$ to $N\rtimes _\alpha G$ with respect to the dual weight (See Lemma 1.4 of Kawahigashi \cite{Kwh}).
Probably, the following proposition is well known for specialists.
\begin{prop}
\label{probability}
In the above context, assume that $G$ is abelian and the action $\alpha $ commutes with the expectation $E$ and $\alpha $ fixes a normal faithful state $\phi_0$ of $N$.
Then the probability index of $\tilde{E}$ is the same as that of $E$.
\end{prop}
\begin{proof}
Since $\tilde{E}$ commutes with the dual action $\hat{\alpha}$, there exists a conditional expectation $\tilde{\tilde{E}}$ from $M\rtimes _\alpha G\rtimes _{\hat{\alpha}}  \hat{G}$ to  $N\rtimes _\alpha G\rtimes  _{\hat{\alpha}} \hat{G}$ with respect to the double dual weight $\hat{\hat{\phi}}$. 
By using Takesaki's duality theorem, it is possible to think that $\tilde{\tilde{E}}$ is an expectation from $M\otimes B(L^2G) \to N\otimes B(L^2G)$. 
Since $\alpha $ fixes the state $\phi :=\phi _0\circ E$, the double dual weight corresponds to $\phi \otimes \mathrm{Tr}$. 
Hence by the uniqueness of the conditional expectation, $\tilde{\tilde{E}}$ corresponds to $E\otimes \mathrm{id}_{B(L^2G)}$. 
Hence we have $\mathrm{Ind}(E)=\mathrm{Ind}(\tilde{E})$.  
\end{proof}

\subsection{Ergodicity of crossed products}
We will classify actions of compact abelian groups. 
In order to achieve this, we need to classify actions of discrete abelian groups on the crossed products. 
The following theorem is useful for investigating the crossed products.
\begin{theo}
\label{ergodicity}
Let $N\subset M$ be an inclusion of semifinite von Neumann algebras with its probability index less than $4$. Let  $A\subset \mathcal{Z}(M)$ be an abelian von Neumann subalgebra.
Identify $A$ with $L^\infty (\Gamma , \mu )$.
Let 
\[ (M \supset N)=\int ^\oplus_\Gamma \left( M(\gamma ) \supset N(\gamma )\right) \ d\mu (\gamma) \]
be the direct integral decomposition of $M\supset N$ over $\Gamma$. 
Let $\theta $ be an action of a locally compact group $G$ on the inclusion $M\supset N$ which is ergodic on $A$. 
Then there exists an inclusion of von Neumann algebras $P\supset Q$ such that $\left( M(\gamma) \supset N(\gamma )\right) \cong (P\supset Q)$ for almost every $\gamma \in \Gamma$. 
\end{theo}
\begin{proof}
Let $T$ be the action of $G$ on $(\Gamma , \mu )$ defined by
\[ \theta _g(f)=f\circ T_{g^{-1}}\]
for $f\in L^\infty (\Gamma , \mu )$, $g\in G$.
Set 
\[ \mathcal{W}:=\{ \mathrm{von} \ \mathrm{Neumann} \ \mathrm{subalgebra} \ \mathrm{of} \ B(H)\}.\]
We consider the Effros Borel structure of $\mathcal{W}$ (For the definition, see Subsection V.6 of Takesaki \cite{T}).
By replacing $H$ by $H\otimes H$, if necessary, we may assume that there exists a subfactor $B$ of $B(H)$ satisfying the following three conditions.

\bigskip

(1) We have $M(\gamma )\subset B$ for almost every $\gamma \in \Gamma$.

(2) The commutant $B'$ is properly infinite.

(3) The factor $B$ is of type $\mathrm{I}_\infty$.

\bigskip

Set
\[ \mathcal{W}_B:=\{ R\in \mathcal{W} \mid R\subset B\}.\]
Then $\mathcal{W}_B$ is the inverse image of the diagonal set of $\mathcal{W}\times \mathcal{W}$ by the Borel map $\mathcal{W} \ni R \mapsto (R\cap B ,R)\in \mathcal{W}\times \mathcal{W}$ (See Corollary IV 8.6 of Takesaki \cite{T}). 
Hence $\mathcal{W}_B$ is Borel. 

\bigskip 

Since $M$ and $N$ are semifinite, there exists a semifinite von Neumann algebra $P$ with $M(\gamma )\cong P$ for almost every $\gamma \in \Gamma$.
Hence by removing a null set, we may assume that $M(\gamma ) \cong P$ for all $\gamma \in \Gamma$ and the map $\gamma \mapsto M(\gamma )$ is Borel measurable. 
Since the commutants of $M(\gamma )$ and $P$ are properly infinite, $M(\gamma )$ and $P$ are spatially isomorphic. 
Hence for almost every $\gamma$, there exists a unitary $u_\gamma $ of $B(H)$ satisfying $u_\gamma M(\gamma )u_\gamma ^*=P$.
We will show that it is possible to choose the map $\gamma \mapsto u_\gamma $ so that it is $\mu$-measurable.

\bigskip 

\textbf{Claim 1}.
There exists a $\mu$-measurable map $\gamma \mapsto u_\gamma \in \mathcal{U}(B(H))$ with $u_\gamma M(\gamma ) u_\gamma =P$ for every $\gamma \in \Gamma$.

\textit{Proof of Claim 1}. 
Set
\[ A:=\{ (\gamma ,u)\in \Gamma \times \mathcal{U}(B(H))\mid u M(\gamma )u^* =P\} .\]
We first show that the set $A$ is Borel.
Since the map $\gamma \mapsto M(\gamma )$ is Borel, by Corollary IV 8.3 of Takesaki \cite{T}, there exists a sequence $\{ f_n :\Gamma \to B(H)\}$ of Borel maps such that for each $\gamma \in \Gamma$, the sequence $\{ f_n(\gamma )\}$ is dense in the unit ball of $M(\gamma )$. 
Then the map $(\gamma , u)\mapsto uf_n(\gamma )u^*$ is Borel and the sequence $\{ uf_n (\gamma )u^*\}$ is dense in $uM(\gamma )u^*$. 
Hence by Corollary IV 8.3 of Takesaki \cite{T}, the map $(u, \gamma ) \mapsto uM(\gamma )u^*$ is Borel.
Hence the set $A$ is Borel.
On the other hand, the map $A\ni (\gamma, u) \mapsto \gamma \in \Gamma $ is surjective. 
Hence by Borel cross section theorem, there exists a $\mu $-measurable map $\gamma \mapsto u_\gamma$ such that $(\gamma, u_\gamma )\in A$ for almost every $\gamma \in \Gamma$.
\qed

\bigskip

Hence by removing a null set, we may assume that the maps $\gamma \mapsto u_\gamma $ and $\gamma \mapsto N(\gamma )$ are Borel. Hence the map $\gamma \mapsto (N(\gamma ), u_\gamma ) \mapsto u_\gamma N(\gamma ) u_\gamma ^*$ (a composition of Borel maps) is also Borel.
Set 
\[ \mathcal{W}_P:=\{ \mathrm{von} \ \mathrm{Neumann} \ \mathrm{subalgebra} \ \mathrm{of} \ P \}.\]
Then $\mathcal{W}_P$ is the inverse image of the diagonal set of $\mathcal{W}\times \mathcal{W}$ by the Borel map $\mathcal{W} \ni Q \mapsto (Q\cap P ,Q)\in \mathcal{W}\times \mathcal{W}$ (See Corollary IV 8.6 of Takesaki \cite{T}). 
Hence $\mathcal{W}_P$ is Borel. 

\bigskip

\textbf{Claim 2}. The map 
\[ \Phi : \mathrm{Aut} (P) \times \mathcal{W}_P \ni (\sigma , Q) \mapsto \sigma (Q) \in \mathcal{W}_P \]
is Borel.

\textit{Proof of Claim 2}. Since the map $\mathcal{W}_P \ni Q \mapsto Q \in \mathcal{W}_P$ is Borel, by Corollary IV 8.3 of Takesaki \cite{T}, there exists a family  of Borel maps $\{ f_n : \mathcal{W}_P \to P_1\subset B(H) _1\}_{n=1}^\infty $ such that for any $Q\in \mathcal{W}_P$, the set $\{ f_n (Q)\mid n=1,2, \cdots \}\subset Q_1$ is weak* dense. 
Here, $P_1$ and $B(H)_1$ are the unit balls of $P$ and $B(H)$, respectively. 
Set
\[ g_n: (\sigma, Q) \mapsto \sigma (f_n (Q))\in \sigma (Q)_1 \subset P_1\subset B(H)_1.\]
Then the map $g_n $ is Borel measurable because the map $\mathrm{Aut}(P) \times P_1 \ni (\sigma , x) \mapsto \sigma (x) \in P_1$ is continuous with respect to the weak* topology. 
On the other hand, for any $(\sigma, Q) \in \mathrm{Aut}(P) \times \mathcal{W}_P $, the set $\{ g_n(\sigma , Q) \mid n=1,2,\cdots \} \subset \sigma (Q)_1$ is weak* dense.
Hence by Corollary IV.8.3 of Takesaki \cite{T}, $\Phi $ is Borel measurable.
\qed

\bigskip

 Now, we  return to the proof of the theorem. 
For $Q\in \mathcal{W}_P$, the set
\[ A_Q:=\{ ( \gamma , \sigma ) \in \Gamma \times \mathrm{Aut}(P) \mid \sigma (u_\gamma N(\gamma )u_\gamma ^{-1})=Q\}\]
is Borel by Claim 2. 
Hence by Corollaries A.8, A.10 of Takesaki \cite{T}, the set
\[ B_Q:=\{ \gamma \in \Gamma \mid \left( N(\gamma ) \subset M(\gamma )\right) \cong (Q \subset P) \} \]
is a Souslin set.
Hence by Theorem A.13 of Takesaki \cite{T}, the set $B_Q$ is $\mu$-measurable.
For any $g\in G$, the set $T_g(B_Q)\bigtriangleup B_Q$ is a null set. 
This is because for almost every $\gamma \in \Gamma$, the inclusion $M(\gamma) \supset N (\gamma ) $ is isomorphic to $ M(T_{g^{-1}}(\gamma )) \supset N(T_{g^{-1}}(\gamma ))$ by the map $x_\gamma \mapsto (\theta _g(x))_{T_{g^{-1}}(\gamma)}$. 
Let $G_0$ be a countable dense subgroup of $G$.
Identify $\mathcal{W}_P$ with $[0,1)$, as  Borel spaces.
By using Ocneanu's classification or Popa's classification, the number of the isomorphism classes of subfactors of type $\mathrm{II}_1$ with index less than $4$ is at most countable (See, Ocneanu \cite{Oc} or Section 5 of Popa \cite{P}).
Hence it is possible to divide $\Gamma$ into at most countably many $B_Q$'s.
Hence by the ergodicity of $T$, there exists $Q\in \mathcal{W}_P$ with $\mu (B_Q)=1$.
\end{proof}

\section{Invariants for actions}
\subsection{Approximate innerness of actions}
\label{approximate innerness}
Let $M\supset N$ be an inclusion of von Neumann algebras. 
Assume that there exists a finite index normal faithful conditional expectation $E$ from $M$ to $N$. 
We investigate when an automorphism $\alpha \in \mathrm{Aut}(M,N)$ is approximated by elements of $\mathrm{Int}(M,N)$ (approximately inner). 
When $M\supset N$ is a strongly amenable inclusion of AFD factors of type $\mathrm{II}_1$,  approximate innerness is characterized by using the Loi invariant (See Theorem 5.4 of Loi \cite{L}).
 Here, we extend his result for inclusions of non-factors (See also Winsl\o w \cite{W}, \cite{W2}).
 Since the expectation $E:M\to N$ has finite index, it is possible to construct the Jones tower $N \subset M \subset M_1 \subset M_2 \subset \cdots $ and the family of Jones projections $\{ e_k\} _{k=1}^\infty$ of the inclusion $M\supset N$ (p.30 of Kosaki \cite{Ko}). 
Assume that the center $\mathcal{Z}(M)$ of $M$ is the same as that of $N$. 
This condition implies that $\mathcal{Z}(M_k)=\mathcal{Z}(M)$ for all $k\geq 1$ because we have
\begin{align*}
\mathcal{Z}(M_1)&=M_1\cap M_1' \\
                &=(J_MN'J_M)\cap (J_MN'J_M)' \\
                &=(J_MN'J_M)\cap (J_MNJ_M) \\
                &=J_M\mathcal{Z}(N)J_M \\
                &=\mathcal{Z}(M). 
\end{align*}
Here, $J_M$ is the modular conjugation.
By the same argument as in Lemma 5.1 of Loi \cite{L} or Lemma 1.4 of Kawahigashi \cite{Kwh}, any automorphism $\alpha $ of $M\supset N$ extends to $M_k$'s by $\alpha (e_k)=e_k$ ($k\geq 1$). 
For an automorphism $\alpha $ of the inclusion $M\supset N$, set
\[\Phi _k(\alpha ):=\alpha |_{M_k\cap N'}\]
for $k\geq 0$. Set
\[ \Phi (\alpha ):= \{ \Phi _k(\alpha )\}_{k=0}^\infty .\]
The sequence $\Phi (\alpha )$ is said to be the Loi invariant of $\alpha$. 
The set 
\[ \mathcal{G}:=\{\Phi (\alpha )\mid \alpha \in \mathrm{Aut}(M,N), \ E\circ \alpha =\alpha \circ E \}\]
has a topological group structure.
The topology is defined by $\Phi (\alpha _n) \to \Phi (\alpha )$ if $\Phi _k(\alpha _n)\to \Phi _k(\alpha)$ for any $k\geq 0$. 
The map $\alpha \mapsto \Phi (\alpha )$ is a group homomorphism. 
\begin{rema}
Let $M \supset N$ be as above. 
Consider the direct integral decomposition 
\[\int ^\oplus _X (M(x) \supset N(x)) \ d\mu (x).\]
If we assume that all fibers $(M(x) \supset N(x)$ are mutually isomorphic to a subfactor $P\supset Q$, then each $\Phi (\alpha ) \in \mathcal{G}$ defines a homomorphism from $X\rtimes _{\alpha}\mathbf{Z}$ to $\Phi (P,Q)$ by the following way.
\[ \Phi (\alpha )_{n,x}:=\Phi (\alpha _{n,x})\]
for $n \in \mathbf{Z}$, $x\in X$. 
\end{rema}
We have the following proposition. 
\begin{prop}
\label{loi}
Let $M \supset N$ be an inclusion of semifinite von Neumann algebras. 
Assume that the inclusion is of the form $(P\supset Q)\otimes \mathcal{Z}(M)$ for a strongly amenable inclusion $P\supset Q$ of semifinite AFD factors. 
Then we have the following two statements.

\textup{(1)} An automorphism $\alpha \in \mathrm{Aut}(M,N)$ is approximately inner if and only if $\Phi (\alpha )=\mathrm{id}$ and the automorphism $\alpha$ admits an invariant trace $\tau$ of $N$.

\textup{(2)} Let $\alpha $, $\beta$ be two actions of a discrete group $G$ on the inclusion $M\supset N$ such that $\Phi (\alpha )$ is conjugate to $\Phi (\beta )$ in $\mathcal{G}$. 
Assume that they have invariant traces $\tau_\alpha $, $\tau _\beta $ of $N$, respectively. 
Then there exists an automorphism $\sigma $ of $M\supset N$ such that $\sigma \circ \alpha _g \circ \sigma ^{-1} \circ \beta _g $ is approximately inner for all $g \in G$.
\end{prop}
\begin{proof}
(1) Identify $\mathcal{Z}(M)$ with $L^\infty (\Gamma, \mu )$. If $\alpha \in \mathrm{Aut}(M,N)$, then the automorphism $\alpha |_{\mathcal{Z}(M)}$ is trivial. Hence we may assume that the automorphism $\alpha |_{\mathcal{Z}(M)}$ is trivial. Let 
\[ M_k=\int ^\oplus _\Gamma M_k(\gamma ) \ d\mu (\gamma ) ,\]
\[ \alpha =\int ^\oplus _\Gamma \alpha _\gamma \ d\mu (\gamma )\]
be the direct integral decomposition.
Then by the same argument as in the proof of Theorem 9.4 of Masuda--Tomatsu \cite{MT}, we have $\alpha \in \overline{\mathrm{Int}}(M,N)$  if and only if $\alpha |_{\mathcal{Z}(M)}$ is trivial and $\alpha _\gamma \in \overline{\mathrm{Int}}(P,Q)$ for almost every $\gamma \in \Gamma$. 
On the other hand, if $\alpha |_{\mathcal{Z}(M)}$ is trivial, then the element $\Phi (\alpha )=\{ \Phi _k(\alpha)\}$ of $\mathcal{G}$ is decomposed into the direct integral.
\[ \Phi (\alpha )=\int ^\oplus _\Gamma \Phi (\alpha ) _\gamma \ d\mu (\gamma )= \int ^\oplus _\Gamma \{(\Phi _k(\alpha ))_\gamma \} \ d\mu (\gamma ).\]
We have $\Phi (\alpha )=\mathrm{id}$ if and only if $\alpha |_{\mathcal{Z}(M)}=\mathrm{id}$ and $\Phi (\alpha ) _\gamma =\mathrm{id}$ for almost every $\gamma \in \Gamma$. 
Hence by Theorem 5.4 of Loi \cite{L}, we have $\Phi (\alpha )=\mathrm{id}$ and $\mathrm{mod}_\tau (\alpha ) =1$ if and only if $\alpha \in \overline{\mathrm{Int}}(M,N)$.

(2) Since the map $\Phi :\mathrm{Aut}(M,N) \to \Phi (\mathrm{Aut}(M,N))$ is surjective, there exists an automorphism $\sigma \in \mathrm{Aut}(M,N)$ with 
\[ \Phi (\sigma \circ \alpha \circ \sigma ^{-1})=\Phi (\beta ).\]
We need to replace $\sigma $ so that it satisfies that $\mathrm{mod}_\tau (\sigma \circ \alpha \circ \sigma ^{-1})=\mathrm{mod}_\tau (\beta)$, where, $\mathrm{mod}_\tau (\alpha)$ is the unique positive operator affiliated with $\mathcal{Z}(M)$ satisfying $\tau \circ \alpha =\tau (\mathrm{mod}_\tau (\alpha )\cdot )$. 
When $M\supset N$ is of type I or finite, then we need not do anything (See lines 11--15 of p.236 of Jones--Takesaki \cite{JT}). 
When $M\supset N$ is of type $\mathrm{II}_\infty$, there exists an inclusion $P_0\supset Q_0$ of finite von Neumann algebras with $\left( (P_0 \supset Q_0 )\otimes R_{0,1}\otimes A\right)=M\supset N$. 
By the same argument as in the proof of Theorem XVIII.2.1 (i) of Takesaki \cite{T3} (p.315), there exists an automorphism $\rho \in \mathrm{Aut}(R_{0,1}\otimes A)$ such that we have
\[ \mathrm{mod}_\tau ((\mathrm{id}_{P_0\supset Q_0}\otimes \rho )\circ \sigma \circ \alpha \circ \sigma ^{-1} \circ (\mathrm{id}_{P_0\supset Q_0}\otimes \rho )^{-1})=\mathrm{mod}_\tau (\beta )\]
and the automorphism $\rho $ is trivial on $A$.
Since $\Phi (\mathrm{id}_{P_0\supset Q_0}\otimes \rho )=\mathrm{id}$, we have the desired conclusion. 
\end{proof}

\begin{rema}
\label{loi2}
Similar results hold when we have one of the following.

(1) $\mathcal{Z}(M) \supset \mathcal {Z}(N)=L^\infty (\Gamma , \mu )$ and $M(\gamma ) \supset N(\gamma ) \cong (\mathbf{C}^p\supset \mathbf{C}) \otimes R$ for some semifinite AFD factor $R$.

(2) $L^\infty (\Gamma , \mu ) \mathcal{Z}(M) \subset \mathcal {Z}(N)$ and $M(\gamma ) \subset N(\gamma ) \cong (M_p(\mathbf{C}) \subset \mathbf{C}^p)\otimes R$ for some semifinite AFD factor $R$.
\end{rema}

\subsection{Centrally trivial part}
In this subsection, we always assume that the inclusion $M\supset N$ is isomorphic to $A \otimes R_{0,1} \otimes (P\supset Q)$. 
Here, $P\supset Q $ is a strongly amenable extremal inclusion of factors of type $\mathrm{II}_1$ satisfying $\overline{\mathrm{Int}}(P,Q)=\mathrm{Aut}(P,Q)$, $A$ is an abelian von Neumann algebra and $R_{0,1}$ is the AFD factor of type $\mathrm{II}_\infty$.

\begin{lemm}
Let $P\supset Q$ be a strongly amenable inclusion of subfactors of type $\mathrm{II}$. Then $\mathrm{Cnt}(P,Q)$ is a Souslin subset of $\mathrm{Aut}(P,Q)$.
\end{lemm}
\begin{proof}
Let $Q \subset P \subset P_1 \subset P_2 \subset \cdots$ be the Jones tower of the inclusion $P \supset Q$ and $\{ x_n \}$ be a dense sequence of the unit ball $S(P)$ of $P$. 
We define a function $F_n^k$ by
\[ F_n^k : (S(P_k) \setminus \{0\}) \times \mathrm{Aut}(P,Q)\ni (a, \alpha ) \mapsto a\alpha (x_n) -x_na\in P_k.\]
Then this map is continuous. 
By Popa's characterization of central triviality (See Theorem 2.14 (6) of Popa \cite{P2}),  $\alpha \in \mathrm{Aut}(P,Q)$ is centrally trivial if and only if there exists $k \in \mathbf{N}$ and $a\in S(P_k)\setminus \{0\}$ with $F_n^k (a, \alpha )=0$ for any $n$. On the other hand, the set 
\[ \bigcap _{n=1}^\infty \{ (a, \alpha ) \in (S(P_k) \setminus \{ 0\} )\times \mathrm{Aut}(P,Q)\mid F_n^k(a, \alpha ) =0\}\]
is a Borel subset. Hence the set
\begin{align*}
A_k:=\{ \alpha \in \mathrm{Aut}(P,Q) & \mid \mathrm{there} \ \mathrm{exists} \ a\in S(P_k ) \setminus \{ 0\} \\
                                     & \mathrm{with} \ F_n^k (a, \alpha )=0 \ \mathrm{for} \ \mathrm{any} \ n\}
\end{align*}
is a Souslin set by Corollary A.10 of Takesaki \cite{T}. Hence $\mathrm{Cnt}(P,Q)=\bigcup _{k=1}^\infty A_k$ is also a Souslin set. 
\end{proof}
Let $G$ be a discrete abelian group and $\alpha $ be a centrally ergodic action of $G$ on $M\supset N$. Set
\[ H(\alpha ):=\{ g\in G\mid \alpha _g |_{A}=\mathrm{id}_A\}.\]
We identify $A$ with $L^\infty (X, \mu )$, and for each $x\in X$, we set
\[N(\alpha , x) := \{ g\in H(\alpha )\mid \alpha _{g,x}\in \mathrm{Cnt}(P,Q)\}.\]
Then the central ergodicity of $\alpha $ and the Souslinness of $\mathrm{Cnt}(P,Q)$ implies that $N(\alpha, x)$ does not depend on the choice of $x\in X$, except for a null set.
Hence we set 
\[ N(\alpha ):= N(\alpha ,x).\]
\begin{theo}
\textup{(Theorem 3.1 of Masuda \cite{M})} \label{masuda} Let $P \supset Q$ be an extremal strongly amenable inclusion of factors of type $\mathrm{II}_1$ with $\mathrm{Aut}(P,Q)=\overline{\mathrm{Int}}(P,Q)$. 
Then for $\sigma \in \mathrm{Cnt}(P,Q)$ and $\theta \in \mathrm{Aut}(P,Q)$, there exists a unique unitary $u(\theta , \sigma ) $ of $N$ with the following two properties.

\textup{(1)} We have $\theta \circ \sigma \circ \theta ^{-1} =\mathrm{Ad}u(\theta, \sigma  ) \sigma$.

\textup{(2)} For any $a\in P_k$ with $\sigma (x)a =ax $ for any $x\in M$, we have $\theta (a)=u(\theta , \sigma )a$.

This $u(\theta , \sigma )$ has the following properties.

\textup{(3)} We have $u(\theta , \sigma _1\circ \sigma _2)=u(\theta , \sigma _1)\sigma _1(u(\theta , \sigma _2))$ for $\theta \in \mathrm{Aut}(P,Q)$, $\sigma _1, \sigma _2\in \mathrm{Cnt}(P,Q)$.

\textup{(4)} We have $u(\theta _1\circ \theta _2 , \sigma )=\theta _1(u(\theta _2 , \sigma ))u(\theta _1, \sigma )$ for $\theta _1$, $\theta _2\in \mathrm{Aut}(P,Q)$, $\sigma \in \mathrm{Cnt}(P,Q)$.
\end{theo}
In our setting, we need to consider automorphisms which may not preserve the center. 
\begin{lemm}
\label{5}
Let $\theta $ be an automorphism of $M \supset N$ which preserves a trace $\tau$. 
Then there exists an automorphism $\overline{\theta} $ of $M\supset N$ which satisfies the following conditions.

\textup{(1)} The automorphism $\overline{\theta}$ preserves the trace $\tau$. 

\textup{(2)} The automorphism $\overline{\theta}$ is of the form $\theta _0\otimes \mathrm{id}_P$, where we identify $M\supset N$ with $(A\otimes R_{0,1})\otimes (P\supset Q)$.
In particular, $\overline{\theta}$ commutes with any automorphism of the form $\mathrm{id}_A\otimes \mathrm{id}_{R_{0,1}} \otimes \alpha$ for some $\alpha \in \mathrm{Aut}(P,Q)$. 

\textup{(3)} We have $\theta \circ \overline{\theta }^{-1}|_A=\mathrm{id}_A$.
\end{lemm}
\begin{proof}
Identify $M\supset N$ with $M\otimes (P\supset Q)$.
Let 
\[ \theta :=\int ^\oplus _X \theta _x \ dx\]
be the direct integral decomposition. 
Set $\overline{\theta }_x:=\theta _x\otimes \mathrm{id}_P$.
Then the automorphism
\[\overline{\theta }:=\int ^\oplus _X \overline{\theta}_x \ dx\]
satisfies conditions (1)--(3).
\end{proof}
\begin{lemm}
\label{6}
Let $\theta $ be a trace-preserving automorphism of $M \supset N$ and $\sigma $ be a centrally trivial automorphism of $P \supset Q$. Set
\[ \sigma := \mathrm{id}_A\otimes \mathrm{id}_{R_{0,1}} \otimes \sigma _0.\]

Then there exists a unique unitary $u(\theta, \sigma )$ of $N$ which satisfies the following two conditions.

\textup{(1)} We have $\theta \circ \sigma \circ \theta ^{-1} =\mathrm{Ad}u(\theta , \sigma ) \sigma $.

\textup{(2)} There exist $k\geq 0$ and $a \in M_k\setminus \{0\}$ with the following properties.
 
\textup{(a)} We have $\sigma _0(x) a=ax$ for any $x\in P$.

\textup{(b)} We have $\theta (a ) =u(\theta , \sigma ) a$.

\textup{(c)} The central support of $a$ is $1$.
\end{lemm}
\begin{proof}
Choose an automorphism $\overline{\theta }$ as in the previous lemma.
Set $\beta :=\theta \circ \overline{\theta}^{-1}$. 
Then $\beta$ is trivial on the center and preserves a trace. 
Hence there exists a unitary $u$ of $N$ such that the automorphism $\alpha :=\mathrm{Ad}u \circ \beta$ preserves $A\otimes \mathbf{C}\otimes (P\supset Q)$ globally.  
Hence by Theorem \ref{masuda}, there exists a unitary $u(\alpha, \sigma )$ of $N$ such that

(1) We have $\alpha \circ \sigma \circ \alpha ^{-1}=\mathrm{Ad}u(\alpha , \sigma )\sigma $.

(2) For any $a\in M_k$ with $\sigma(x)a=ax$ for any $x\in M$, we have $\alpha (a)=u(\alpha , \sigma )a$. 

Set $u(\theta, \sigma ):=u^*u(\alpha, \sigma )\sigma (u)$. We will show that this $u(\theta, \sigma)$ satisfies the desired conditions. 
Condition (1) is shown by 
\begin{align*}
\theta \circ \sigma \circ \theta ^{-1} &= \mathrm{Ad}u^* \circ \alpha \circ \overline{\theta } \circ \sigma \circ \overline{\theta }^{-1} \circ \alpha ^{-1}\circ \mathrm{Ad}u \\
                                       &=\mathrm{Ad}u^*\circ \alpha \circ \sigma \circ \alpha ^{-1} \circ \mathrm{Ad}u \\
                                       &=\mathrm{Ad}(u^*u(\alpha , \sigma )\sigma (u) )\circ \sigma .
\end{align*}
Let $a_0\in P_; \setminus \{ 0\}$ be an element satisfying $\sigma _0 (x)a_0=a_0x$ for any $x\in P$.
Set $a=1\otimes 1\otimes a_0\in A\otimes B(l^2) \otimes P=M$.
Then for any $x\in M$, we have $\sigma (x)a=ax$ and the central support of $a$ is $1$.
We need to show condition (b). 
This is shown by
\begin{align*}
\theta (a)&=\theta \circ \overline{\theta }^{-1} (a) \\
                             &=\mathrm{Ad}u^*\alpha (a) \\
                             &=u^*u(\alpha , \sigma )au \\
                             &=u^*u(\alpha , \sigma )\sigma (u)a.
\end{align*}

Next, we show the uniqueness. Assume that another unitary $v(\theta , \sigma )$ satisfies the above two conditions. 
Then by the first condition, there exists a unitary $c$ of $A$ with $v(\theta, \sigma ) =cu(\theta, \sigma )$. 
However, by the second condition, this $c$ should be $1$.  
\end{proof}

For an automorphism $\rho $ of $M\supset N$ which is trivial on $A$, it is possible to construct a unitary $u(\rho , \sigma )$ of $N$ in the following way. 
Let
\[ \rho =\int ^\oplus _X \rho _x \ dx\]
be the direct integral decomposition.
For each $x\in X$, we identify $(M_x\supset N_x)\cong  R_{0,1} \otimes (P\supset Q)$ with $ M_x\otimes (M_x\supset N_x)$. Set $\overline{\rho}_x:=\rho _x\otimes \mathrm{id}_{M_x}$.
Set 
\[\overline {\rho}:=\int ^\oplus _X \overline{\rho}_x \ dx.\]
Then the automorphism $\rho \circ \overline{\rho}^{-1}$ is trace-preserving. 
Hence by using the above lemma, we obtain a unitary $u(\rho \circ \overline{\rho}^{-1}, \sigma )$. 
Set $u(\rho , \sigma ):=u(\rho \circ \overline{\rho}^{-1}, \sigma )$. 
Since $\overline{\rho}_x^{-1}=\rho ^{-1}_x \otimes \mathrm{id}_{M_x}$ commutes with $\sigma _x$, this $u(\rho , \sigma )$ satisfies the following conditions.

\bigskip

(1) We have $\rho \circ \sigma \circ \rho ^{-1}=\mathrm{Ad}u(\rho , \sigma ) \sigma $.

(2) For any $a \in M_k$ such that $\sigma (x)a=ax$ for any $x\in M$, we have $\rho (a)=u(\rho, \sigma )a$.

\bigskip

By using this unitary, it is possible to show the following result.

\begin{cor}
\label{coro}
Let $\theta $ be a trace-preserving automorphism of $M\supset N$ and $\sigma _0$ be a centrally trivial automorphism of $ P \supset Q$. 
Then for any automorphism $\rho  $ of $M\supset N $ which is trivial on $A$, we have
\[ u(\rho \circ \theta \circ \rho ^{-1} , \sigma )=\rho (\theta (u(\rho ^{-1}, \sigma ))u(\theta , \sigma ))u(\rho , \sigma ).\] 
\end{cor}
\begin{proof}
This corollary is shown by the uniqueness of $u(\theta, \sigma )$.
\end{proof}

\subsection{Characteristic invariants}
In this subsection, we always assume that the inclusion $M \supset N$ satisfies the same assumption as that of the previous subsection.

Set 
\[ \chi (P,Q):=\frac{\overline{\mathrm{Int}}\cap \mathrm{Cnt}}{\mathrm{Int}}(P,Q),\]
\[\chi (M,N):=\frac{\overline{\mathrm{Int}}\cap \mathrm{Cnt}}{\mathrm{Int}}(M,N).\]
If we identify $A$ with $L^\infty(X, \mu)$, it is possible to think of $\chi (M,N)$ as the set of all $\mu $-measurable functions from $X$ to $\chi (P,Q)$. 
In this subsection, we also put the following assumption.

\bigskip

\textbf{Assumption.} 

(1) The group $\chi (P,Q)$ is a finite abelian group. 

(2) There exists a unique right inverse $\sigma ^0:\chi (P,Q)\to \mathrm{Cnt}(P,Q)$ of the quotient, up to cocycle conjugacy.

\bigskip

\begin{rema}
(1) A right inverse $\sigma ^0$ is outer, as an action on $Q$. 
Hence it is minimal because $\chi (P,Q)$ is a finite group. 
Hence the uniqueness up to cocycle conjugacy implies the uniqueness of up to conjugacy.

(2) Under this assumption, any action $\sigma :\chi (P,Q) \to \mathrm{Cnt}(M,N)$ with the following $\mathrm{(*)}$ satisfies the assumption in Lemma \ref{6}.

\bigskip

$\mathrm{(*)}$ For any $s\in \chi (P,Q)$, we have $[\sigma _s]=s$ (a constant function from $X$ to $\chi (P,Q)$).

\bigskip

(3) Some important subfactors satisfy the above assumption. 
For example, Jones subfactors with principal graph $A_{2n+1}$, $n\geq 2$ and subfactors with principal graph $E_6$, $E_8$ satisfy this assumption 
(See Examples 6.1, 6.2, 6.3 and Theorem 6.5  of Kawahigashi \cite{Kwh3}). 
\end{rema}

\bigskip

\begin{lemm}
Let $\alpha $ be a centrally ergodic action of a discrete abelian group $G$. 
Then for any $h\in N(\alpha)$, the element $[\alpha _h]\in \chi (M,N)$ is a constant function, if we think of $[\alpha _h]$ as a function from $X$ to $\chi (P,Q)$.
\end{lemm}
\begin{proof}
By assumption, we have $\mathrm{Aut}(P,Q)=\overline{\mathrm{Int}}(P,Q)$. Hence for any $g\in G$, we have
\begin{align*}
[\alpha _{h,x}] &=[\alpha _{g, x} \circ \alpha _{h,x} \circ (\alpha _{g,x})^{-1}] \\
                &=[\alpha _{g,x} \circ \alpha _{h,x} \circ \alpha _{g^{-1}, g.x}] \\
                &=[\alpha _{ghg^{-1}, g.x}] \\
                &=[\alpha _{h,g.x}].
\end{align*}
In the first equation, we used Theorem \ref{masuda}. This computation shows that the map $x\mapsto [\alpha _{h,x}]$ is $\alpha |_A$-invariant. Thus by the ergodicity of $\alpha |_A$ and the Souslinness of $\mathrm{Cnt}(P,Q)$, the map is a constant function.
\end{proof}
Set
\[ \nu _\alpha :N(\alpha ) \ni h\mapsto [\alpha _h]\in \chi (M,N).\]
By the above lemma, it is possible to think that the range of the map $\nu _\alpha $ is contained in $\chi (P,Q)$. For each $h\in N(\alpha)$, set
\[ \sigma _h:=\mathrm{id}_{A\otimes R_{0,1}}\otimes \sigma ^0_{\nu _\alpha(h)}.\] 
Then for each $h\in N(\alpha )$, there exists a unitary $u^\alpha _h \in N$ with $\alpha _h=\mathrm{Ad}u^\alpha _h\sigma _{\nu _\alpha (h)}$.
\begin{lemm}
For each $g\in G$, $h\in N(\alpha)$, we have
\[ \lambda ^\alpha (g,h):=(u_h^\alpha )^*\alpha _g(u^\alpha _h)u(\alpha_g, \sigma _h) \in A.\]
\end{lemm}
\begin{proof}
This is shown by the equation
\[ \mathrm{Ad}u^\alpha _h \circ \sigma _{\nu _\alpha (h )}=\alpha _h =\alpha _g\circ \alpha _h\circ \alpha _{g^{-1}}=\mathrm{Ad}\alpha _g(u^\alpha _h)u(\alpha _g, \sigma _h)\circ \sigma _{\nu _\alpha (h)}.\]
\end{proof}
\begin{lemm}
For each $n,m \in N(\alpha)$, we have
\[ \mu ^\alpha (n,m):=u^\alpha _n \sigma _{\nu _{\alpha}(m)}(u_n)u_{nm}^*\in A.\]
\end{lemm}
\begin{proof}
This is shown by $\alpha _n\alpha_m =\alpha _{nm}$.
\end{proof}
\begin{prop}
The above $\lambda ^\alpha $ and $\mu ^\alpha $ satisfy the following equations.

\textup{(1)} $\mu ^\alpha (h,k)\mu ^\alpha(hk,l) =\mu ^\alpha (k,l)\mu ^\alpha (h,kl)$ for $h,k,l \in N(\alpha)$.

\textup{(2)} $\lambda^\alpha (g_1g_2,n)=\lambda^\alpha (g_1,n) \lambda ^\alpha (g_2, n)$ for $g_1, g_2\in G$, $n\in N(\alpha)$.

\textup{(3)} $\lambda ^\alpha (h,k)=\mu ^\alpha (h,k)\mu ^\alpha (k,h)^*u(\sigma _{\nu _\alpha(k)}, \sigma _{\nu _\alpha (h)})^*$ for $k,h\in N(\alpha )$.

\textup{(4)} $\lambda ^\alpha (g, hk)\lambda ^\alpha (g,h)^*\lambda ^\alpha (g,k)^*=\mu ^\alpha (h,k)\mu ^\alpha (h,k)^*$ for $g\in G$, $h,k\in N(\alpha)$.

\textup{(5)} $\lambda ^\alpha (g,1)=\lambda ^\alpha (1,h)=\mu ^\alpha (h,1)=\mu ^\alpha (1,h)=1$ for $g\in G$, $h\in N(\alpha)$.
\end{prop}
\begin{proof}
This is shown by the same method as that of Section 4 of Masuda \cite{M}.
\end{proof}
Set
\begin{align*}
 Z(G, N(\alpha ), \kappa _\alpha ):= \{ (\lambda , \mu ) \mid & \lambda :G\times N(\alpha ) \to \mathcal{U}(A), \\ 
                                                              & \mu :N(\alpha )\times N(\alpha ) \to \mathcal{U}(A) \ \mathrm{satisfying} \ \mathrm{the} \ \mathrm{above} \ \mathrm{(1)-(5)}\}.
\end{align*}

We introduce an equivalence relation on $Z(G, N(\alpha ), \kappa _\alpha )$ by the following way.
For $(\lambda ^1, \mu ^1)$, $(\lambda ^2, \mu ^2)\in Z(G, N(\alpha ), \kappa _\alpha )$, we have
\[ (\lambda ^1, \mu ^1) \sim (\lambda ^2, \mu ^2)\]
if  there exists a function $c:N(\alpha ) \to \mathcal{U}(A)$ with $c_1=1$ satisfying the following two equations.

(1) $c_n^*c_m^*\mu ^2(m,n)c_{mn}=\mu ^1(m,n)$ for $m,n \in N(\alpha )$. 

(2) $c_n^*\lambda ^2(g,n)\alpha _g(c_n)=\lambda ^1(g,n)$ for $g\in G$, $n\in N(\alpha )$.

\bigskip

Set $\Lambda (G, N(\alpha )|\kappa _\alpha ):=Z(G, N(\alpha ), \kappa _\alpha )/\sim$. 

\begin{prop}
\label{invariant}
Let $\alpha $, $\beta$ be two actions of a discrete abelian group $G$ on $M\supset N$ which preserve traces $\tau _\alpha $, $\tau _\beta $, respectively.
 Assume that there exists an automorphism $\theta $ of $M\supset N$ and an $\alpha $-cocycle $\{u_g\} \subset N$ satisfying $\mathrm{Ad} u_g \circ \alpha =\theta \circ \beta \circ \theta ^{-1}$. 
Then we have $[\theta (\lambda ^\beta, \mu ^\beta)]=[(\lambda ^\alpha, \mu ^\alpha)]$ in $\Lambda (G, N(\alpha )|\kappa _\alpha)$.
\end{prop}
\begin{proof}
Since $\alpha $ and $\beta $ are mutually cocycle conjugate, we have $N(\alpha )=N(\beta)$.
 By the same argument as that just after Lemma \ref{6}, there exists an automorphism $\overline{\theta}$ of $M\supset N$ which satisfies the following conditions.
 
 \bigskip

(1) The automorphism $\overline{\theta}$ preserves $\tau_\beta$.

(2) We have $\theta |_A=\overline{\theta}|_A$.

(3) The automorphism $\overline{\theta}$ commutes with $\sigma _{\nu _\alpha (h)}$  for any $h\in N(\alpha )$.

\bigskip

Then for each $h\in N(\alpha )$, we have $\overline{\theta} \circ \beta_h \circ \overline{\theta}^{-1}=\mathrm{Ad}\overline{\theta}(u^\beta_h) \sigma _{\nu _\alpha(h)}$. 
On the other hand, it is possible to choose $a\in M_k$ with the following properties.

\bigskip

(1) We have $\overline{\theta}(a)=a$.

(2) For any $x\in M$, we have $\sigma (x)a=ax$.

(3) The central support of $a$ is $1$.

\bigskip

Then by the uniqueness of $u(\overline{\theta}\circ \beta _g\circ \overline{\theta}^{-1}, \sigma _{\nu _\alpha (h)}) )$, we have
\[ u(\overline{\theta}\circ \beta _g \circ \overline{\theta}^{-1}, \sigma _{\nu _\alpha (h)})=\overline{\theta}(u(\beta _g, \sigma _{\nu _\alpha (h)})).\]
Hence  we have $[(\lambda ^{\overline{\theta}\circ \beta \circ \overline{\theta}^{-1}}, \mu ^{\overline{\theta}\circ \beta \circ \overline{\theta}^{-1}})]=[\theta (\lambda ^\beta, \mu ^\beta)]$. 
Hence we may assume that $\alpha |_A=\beta |_A$ and there exists an automorphism $\rho $ of $M\subset N $ which is trivial on $A$ and satisfies $\mathrm{Ad} u_g \circ \alpha =\rho \circ \beta \circ \rho ^{-1}$.
Hence by the same argument as that of the proof of Proposition 4.2 and Remark just after Proposition 4.2 of Masuda \cite{M}, we get the conclusion (in the computation of $\lambda ^{\rho \circ \beta \circ \rho ^{-1}}$, we use Corollary \ref{coro}). 
\end{proof}

\begin{rema}
\label{other}
We will classify actions of compact abelian groups on subfactors with index less than $4$, and our strategy is to classify dual actions. 
However, some crossed products by actions of compact abelian groups are not of the above form (for the proof, see Section \ref{compact}). 
We need to handle inclusions of the form $A\otimes R\otimes (P\supset Q)$, where $A$ is abelian, $R$ is a semifinite factor and $P\supset Q$ is an inclusion with one of the following properties.

(1) The inclusion $P\supset Q $ is a subfactor of type $\mathrm{II}_1$ which satisfies $\mathrm{Aut}(P,Q)=\overline{\mathrm{Int}}(P,Q)$ and there exists a unique right inverse $\sigma^0:\chi (P,Q)\to \mathrm{Cnt}(P,Q)$.

(2) The inclusion $P\supset Q $ is a strongly amenable, extremal subfactor of type $\mathrm{II}_1$ which satisfies $\mathrm{Cnt}(P,Q)=\mathrm{Int}(P,Q)$.

(3) The inclusion $P\supset Q $ is a strongly amenable, extremal subfactor of type $\mathrm{II}_1$ with its principal graph of type $D_4$.

(4) $Q$ is a finite factor and $P=Q\otimes \mathbf{C}^p$, where $p=2$ or $3$.

(5) $Q=R\otimes \mathbf{C}^p$ where $R$ is a finite factor and $p=2$ or $3$, $P=R\otimes M_p(\mathbf{C})$.

In case (1) and $R$ is of type $\mathrm{II}_1$, invariants are constructed by the same way. A proposition corresponding to Proposition \ref{invariant} is also shown by the same way. 
In  cases (2)--(5), any centrally trivial automorphism of $P\supset Q$ is implemented by a unitary of $P$ or $P_1$, where $P_1$ is the basic extension of $P\supset Q$. 
Hence by using these unitaries, it is possible to construct characteristic invariants by the same argument as that of the constructions in Jones--Takesaki \cite{JT}.
\end{rema}

\section{Classification of actions of discrete abelian groups}
In this section, we always assume that $M \supset N$ is an inclusion of the form $A\otimes B(l^2) \otimes (P\supset Q)$, where $A$ is abelian, $B(l^2)$ is a type I factor and $P\supset Q$ is an inclusion of factors of type $\mathrm{II}_1$ which is strongly amenable and extremal. 
We also assume that the inclusion $P\supset Q$ satisfies one of the following conditions.

(1) Any automorphism of $P\supset Q$ is approximately inner and there exists a unique right inverse from $\chi (P,Q)$ to $\mathrm{Cnt}(P,Q)$, up to (cocycle) conjugacy.

(2) Any centrally trivial automorphism is inner.

\begin{theo}
\label{classify}
Let $\alpha $ and $\beta $ be two centrally ergodic actions of a discrete abelian group $G$ on $M \supset N$. 
Assume that they preserve traces, respectively. 
Then they are mutually cocycle conjugate if and only if there exists an automorphism $\theta $ of $M\supset N$ with $\Phi (\theta \circ \alpha \circ \theta ^{-1}) =\Phi (\beta)$, $N(\alpha )=N(\beta )$, $\nu _\alpha =\nu _\beta$ and $[\theta (\lambda ^\alpha , \mu ^\alpha )]=[(\lambda ^\beta , \mu ^\beta)]$.
\end{theo}
The ``only if'' part has been shown in the end of the previous section. 
Hence we will show the ``if'' part. 
As we will explain in Remark \ref{end}, the second case can be shown by the same argument as that  of the proof of Theorem 2.1.14 of Jones--Takesaki \cite{JT}. 
Here we show the theorem when the inclusion $P\supset Q $ satisfies the first condition.
The outline of the proof is the same as that of the case  of single factors (Theorem 2.1.14 of Jones--Takesaki \cite{JT} and Theorem 2.2 of Kawahigashi--Takesaki \cite{KT}). 
In order to imitate their proof, the basic strategy is to use analogies between automorphisms of single factors of type $\mathrm{III} $ and those of subfactors of type $\mathrm{II}_1$ (See Section 3 of Kawahigashi \cite{Kwh2} or Section 3 of Masuda \cite{M}). 
In particular, we use the right inverse $\sigma ^0 :\chi (P,Q) \to \mathrm{Cnt}(P,Q)$ instead of modular automorphism group $\{ \sigma ^\phi _t\}_{t\in \mathbf{R}}$. 
Although the outline of the proof is similar to that of the case of single factors, some parts are technically different.
In particular, the proof of Lemma \ref{point}, which corresponds to Lemma 2.5.6 of Jones--Takesaki \cite{JT}, Lemma 2.4 of Kawahigashi--Takesaki \cite{KT}, is difficult.
Hence we explain how to deal with these problems.

\bigskip

 Let $\alpha $ and $\beta $ be actions of a discrete abelian group $G$ on $M \supset N$ satisfying assumptions in Theorem \ref{classify}. 
Assume that there exists an automorphism $\theta $ of $M\supset N$ with $\Phi (\theta \circ \alpha \circ \theta ^{-1}) =\Phi (\beta)$, $N(\alpha )=N(\beta )$, $\nu _\alpha =\nu _\beta$ and $[\theta (\lambda ^\alpha , \mu ^\alpha )]=[(\lambda ^\beta , \mu ^\beta)]$.
Identify $B(l^2)\otimes (P\supset Q)$ with $(B(l^2)\otimes P)\otimes (B(l^2)\otimes (P\supset Q))$.
 By considering
\[ \int ^\oplus _X \theta ^{-1}_{T, T^{-1}x}\otimes \theta _{T,T^{-1}x} \ dx\]
instead of considering $\theta$, we may assume that $\theta $ preserves the trace which is fixed by $\alpha $ (Here, $T$ is a transformation of $X$ induced by $\theta |A$). 
Then $\theta \circ \alpha \circ \theta ^{-1}$ preserves a trace. 
Hence by considering $\theta \circ \alpha \circ \theta ^{-1}$ instead of considering $\alpha$, we may assume that $\alpha $ and $\beta $ have the same invariants.
 By the same argument as that of the proof of Theorem XVIII.2.1 (i) of Takesaki \cite{T3} (p.315), there exists an automorphism $\rho $ of $M\supset N$ such that $\rho \circ \alpha \circ \rho ^{-1}$ and $\beta $ preserve the same trace. 
Hence we may assume that $\alpha _g\circ \beta _{g^{-1}}$ is approximately inner for any $g\in G$.

Set 
\[ H(\alpha ):= \{ g\in G \mid \alpha _g|_{A}=\mathrm{id}\},\]
\[ K:=G/H(\alpha ).\]
Then $\alpha $ defines an action $\overline{\alpha }$ of $K$ on $A\cong L^\infty(X, \mu )$. Set
\[\mathcal{G}:=X\rtimes _\alpha G, \ \mathcal{K}:=X\rtimes _{\overline{\alpha}}K.\]
Then by Lemma 2.2.11 of Jones--Takesaki \cite{JT}, we have
\[ \mathcal{G}\cong H(\alpha ) \times \mathcal{K}\]
as groupoids. 
The action $\alpha $ of $G$ on $M\supset N$ defines an action $\alpha $ of $\mathcal{G}$ on $(P\supset Q) \otimes B(l^2)$ by the following.
\[ \alpha _{g,x}(a_x)=(\alpha _g(a))_{g.x}\]
for $a=\int ^\oplus _Xa_x dx$, $(g,x) \in \mathcal{G}$.

\begin{lemm}
\label{Sutherland--Takesaki}
\textup{(Theorem 3.1 (b) of Sutherland--Takesaki \cite{ST})}. 
There exists an action $l$ of $H(\alpha ) \times \mathcal{K}$ on the AFD factor $R_{0}$ of type $\mathrm{II}_1$ with $\chi _l =1\in \Lambda (H(\alpha ) \times \mathcal{K}, N(\alpha ), \mathbf{T})$ such that for each homomorphism $q:\mathcal{K}\to \hat{H}$, there exist a Borel map $\rho :X\ni x\mapsto \rho _x\in \mathrm{Aut}(R_0)$ and Borel cocycles $\{ u_k\}_{k\in \mathcal{K}}$, $\{ v_{h,x} \}_{(x,h)\in X\times H(\alpha)}$ with the following conditions.

\textup{(1)} We have $\rho _{r(\gamma )} \circ l_\gamma \circ \rho _{s(\gamma )}^{-1}=\mathrm{Ad}v_\gamma \circ l_\gamma $ for $\gamma \in \mathcal{G}$.

\textup{(2)} We have $u_k l_{1,k}(v_{h, s(k)})=\langle h, q(k)\rangle v_{h, r(k)}l_{h, r(k)}(u_k)$ for $k \in \mathcal{K} $, $h\in H(\alpha )$.
\end{lemm}

For each $(g,x) \in \mathcal{G}$, set
\[ m_{g,x}:=\alpha _{g,x} \otimes l_{g,x} \in \mathrm{Aut}((P\supset Q) \otimes B(l^2) \otimes R_0).\]
Then $m$ defines an action $m$ of $G$ on $M\supset N$.

\begin{lemm}
Let $\alpha $ and $m$ be as above. 
Then we have $\Phi (m)=\Phi (\alpha )$, $N(m)=N(\alpha )$, $\nu _m =\nu_\alpha $ and $[(\lambda ^m, \mu ^m)]=[(\lambda ^\alpha , \mu ^\alpha )]$.
\end{lemm}
\begin{proof}
The only non-trivial assertion in the lemma is that we have $[(\lambda ^m, \mu ^m)]=[(\lambda ^\alpha , \mu ^\alpha )]$. 
By the uniqueness of $\sigma ^0 $, we may assume that $\sigma $ is of the form $\overline{\sigma } \otimes \mathrm{id}_{R_0}$, where $\overline{\sigma }=\mathrm{id}_{A\otimes B(l^2)} \otimes \sigma ^0 $ (See also Proposition \ref{invariant}).
Then by the uniqueness in Theorem \ref{masuda}, we have 

\[ u(\alpha _{g,g^{-1}.x}\otimes \mathrm{id}_{R_0}, \overline{\sigma }_{\nu _\alpha (h), g^{-1}.x} \otimes \mathrm{id}_{R_{0,1}})=u(\alpha_{g,g^{-1}.x},\overline{\sigma }_{\nu _\alpha (h),g^{-1}.x})\otimes 1,\]
\[ u(\mathrm{id}\otimes l_{g,g^{-1}.x}, \overline{\sigma }_{\nu _\alpha (h), g^{-1}.x}\otimes \mathrm{id}_{R_{0,1}})=1.\]
Hence by the uniqueness in Theorem \ref{masuda}, we have
\[ u(m_g , \sigma )=\int ^\oplus _Xu(\alpha _{g,g^{-1}.x}, \overline{\sigma }_{\nu _\alpha (h), g^{-1}.x} )\otimes 1 \ dx.\]
For each $h\in N(\alpha )$, it is possible to choose $u^m_h$ so that
\[ u^m_h=\int ^\oplus _X u^\alpha _{h,x}\otimes u^l_{h,x} \ dx.\]
Hence we have
\[ \lambda ^m(g,h)=\int ^\oplus _X \lambda ^\alpha (g,h)_x \otimes \lambda ^l(g,h)_x \ dx,\]
\[\mu ^m (n,m)=\int ^\oplus _X \mu ^\alpha (n,m)_x \otimes \mu ^l(n,m)_x \ dx.\]
Since $\chi _l=1$, we get the conclusion.
\end{proof}

Hence it is enough to show that the action $m$ is cocycle conjugate to the action $\beta$, and in the rest of this section, we will show this cocycle conjugacy.

\begin{lemm}
\label{factors}
The action $m|_{H(\alpha)}$ is cocycle conjugate to the action $\beta |_{H(\alpha)}$.
\end{lemm}
\begin{proof}
This follows from a classification theorem of actions on subfactors (Theorems 5.1 and 6.1 of Masuda \cite{M} or Section 8 of Masuda \cite{M2}).
\end{proof}
By this lemma, we may assume that there exists an $m$-cocycle $\{v_h\}_{h\in H(\alpha )}$ of $N$ such that $\beta _h=\mathrm{Ad}v_h \circ m_h$ for $h \in H(\alpha )$.

Let 
\[ v_h =\int ^\oplus _X v_{h,x} \ dx\]
be the direct integral decomposition.
Set
\begin{align*}
 \mathcal{A}_0(x):=\{ &(\rho ^x, w^x_h) \mid \rho ^x\in \mathrm{Aut} (M_x, N_x ), \\
                      & \{ w_h^x\} _{h\in H(\alpha )} \ \mathrm{is} \  \mathrm{an} \ m|_{H(\alpha)} \ \mathrm{cocycle}, \ \rho ^x \circ m_{h,x} \circ (\rho ^x)^{-1}=\mathrm{Ad}w^x_h \circ m_{h,x}\},
 \end{align*}
\[ \mathcal{B}_0(x):=\{ (\mathrm{Ad}u^x, u^xm_{h,x}(u^x)^*) \mid u^x \in \mathcal{U}(N_x) \}.\]
Then as in subsection 2.5 of Jones--Takesaki \cite{JT}, it is possible to introduce a Polish group structure on $\mathcal{A}_0(x)$.
The product is defined by the following.
\[ (\rho _1^x, w_h^x)(\rho ^x_2, v_h^x):=(\rho_1^x\circ \rho _2^x , \rho ^x_1(v_h^x)w_h^x).\]
For the meaning of this group structure, see Section 2 of Kawahigashi--Takesaki \cite{KT}.
A sequence $\{ (\rho ^x_n, w^x_{h,n})\}$ converges to $(\rho ^x, w_h^x)$ if, by definition, $\rho ^x_n \to \rho ^x$ and $w^x_{h,n} \to w_h^x$ for any $h\in H$.

Since $H(\alpha )$ is a commutative group, it is possible to think of $\hat{H}$ as a subgroup of $\mathcal{A}_0(x)$ by the following way.
\[ \hat{H}\ni p\mapsto (\mathrm{id}, p) \in \mathcal{A}_0(x).\]
Set
\[ \mathcal{A}(x):=\mathcal {A}_0(x)/\hat{H}.\] 
Let $\pi : \mathcal{A}_0(x) \to \mathcal{A}(x)$ be the quotient. For each $k\in \mathcal{K}$, $k:x\mapsto y$, set
\[ n_k:=\beta_k\circ m_k ^{-1}.\]
Since $\beta _{h,x} =\mathrm{Ad}(v_{h,x})m_{h,x}$, by the same computation as that preceding Lemma 2.4 of Kawahigashi--Takesaki \cite{KT}, we have
\[ n_k\circ m_{h,x}\circ n_k^{-1}=\mathrm{Ad}(\beta _k (v_{h,x})^* v_{h,y})m_{h,y}.\]
Set
\[ w_h^k :=\beta _k ((v_{h,x})^*)v_{h,y}.\]
Then the map $h\mapsto w_h^k$ is an $m_{h,y}$-cocycle. 
By construction, we have $[(n_k, w_h^k)]\in \mathcal{A}(y)$. 
We will show the following lemma, which corresponds to Lemma 2.5.6 of Jones--Takesaki \cite{JT} or Lemma 2.4 of Kawahigashi --Takesaki \cite{KT}.

\begin{lemm}
\label{point}
We have $[(n_k, w_h^k)]\in \overline{\mathcal{B}(y)}$.
\end{lemm}

\begin{lemm}
Lemma \ref{point} implies Theorem \ref{classify}.
\end{lemm}
\begin{proof}
This is shown by the same argument as that of the proof of Theorem 2.1.14 of Jones--Takesaki \cite{JT} or Theorem 2.2 of Kawahigashi-Takesaki \cite{KT}. Lemma \ref{Sutherland--Takesaki} is used here.
\end{proof}
Hence in the rest of this section, we will show Lemma \ref{point}.
\begin{lemm}
\label{inner}
Let $[(n, w_h)]$ be an element of $\mathcal{A}(y)$ and $u $ be a unitary of $N_y$. 
Then $[(n, w_h)] \in \overline{\mathcal{B}(y)}$ if and only if $[(\mathrm{Ad}u\circ n, uw_hm_{h,y}(u)^*)]\in \overline{\mathcal{B}(y)}$. 
\end{lemm}
\begin{proof}
This is trivial.
\end{proof}
By replacing $(n_k, w_h^k)$ by $(\mathrm{Ad}v\circ n_k, vw_h^km_{h,y}(v^*))$ for some unitary $v$ of $N_y$ (and by the above lemma), we may assume that the automorphism $n_k$ is of the from $n'_k\otimes \mathrm{id}_{B(l^2)}$, where we identify $M_y\supset N_y$ with $(P\supset Q)\otimes B(l^2)$.
By Theorem 6.1 of Masuda \cite{M},  the automorphism $n_k$ is cocycle conjugate to $ n_k \otimes \mathrm{id}_{R_0}$, where $R_0 $ is the AFD factor of type $\mathrm{II}_1$. 
Hence by considering $(\mathrm{Ad}u \circ n_k, uw_h^km_{h,y}(u)^*)$ for some unitary $u$ of $N_y$ (and by the above lemma), we may assume that $n_k$ is of the form $\theta ' \circ (n'_k \otimes \mathrm{id}_{R_0}\otimes \mathrm{id}_{B(l^2)})\circ (\theta ') ^{-1}$ for some automorphism $\theta '$ of $M_y\supset N_y$.
As we have explained, in order to show the classification theorem, we will consider a subfactor of type $\mathrm{II}_1$ as a similar object to a single factor of type III. 
In this viewpoint, the right inverse $\sigma ^0 :\chi (P,Q) \to \mathrm{Cnt}(P,Q)$ corresponds to  a modular automorphism group. 
In the case of factors of type III, any automorphism extends to the crossed product by a modular automorphism group, which is said to be the canonical extension (See Section 12 of Haagerup--St\o rmer \cite{HS}). 
In the case of subfactors, any approximately inner automorphism extends to the crossed product by $\sigma ^0$ (See Proposition 3.2 of Masuda \cite{M}). 
We will use this extension.

\begin{prop}
\textup{(Proposition 3.2 of Masuda \cite{M})}. 
\label{3.2 of masuda}
Let $P\supset Q$ be a strongly amenable, extremal inclusion of factors of type $\mathrm{II}_1$ such that any automorphism of $P\supset Q$ has trivial the Loi invariant.
 Assume that there exists a right inverse $\sigma ^0:\chi (P,Q) \to \mathrm{Cnt}(P,Q)$. 
Then any approximately inner automorphism $n$ of $P\supset Q$ extends to an automorphism $\tilde{n}$ of the crossed product  $(P\supset Q) \rtimes _{\sigma ^0 }\chi (P,Q)$ by $\lambda _s\to u(n, \sigma _s)\lambda _s$ for $s\in \chi (P,Q)$. 
The extension $\alpha \mapsto \tilde{\alpha }$ has the following properties. 

\textup{(1)} For any unitary $u$ of $Q$, we have $\tilde{\mathrm{Ad}}u=\mathrm{Ad}u$.

\textup{(2)} The map $\alpha \mapsto \tilde{\alpha }$ is a continuous group homomorphism.
\end{prop}
\begin{proof}
All the statements except for property (2) are shown in Masuda \cite{M}. Here, we show property (2). 
By Theorem \ref{masuda}, for any $a\in P_k$ with $\sigma (x)a=ax$, the map $\alpha \mapsto \alpha (a)a^*= u(\alpha, \sigma )aa^*$ is continuous. 
Hence the map $\alpha \mapsto E_k(u(\alpha , \sigma )aa^*)=u(\alpha , \sigma )E_k(aa^*)$ is continuous, where $E_k:P_k\to P$ is the conditional expectation.
 Thus we are done.  
\end{proof}

We apply the above result for our setting.
Let $\sigma :\chi (P,Q)\to \mathrm{Cnt}(M,N)$ be the action with $[\sigma _s]=s$ for any $s\in \chi (P,Q)$.
Then by assumption, we may assume that the action $\sigma $ is of the form $\sigma ^0\otimes \mathrm{id}_{B(l^2)\otimes A}$, where $\sigma ^0 :\chi (P,Q)\to \mathrm{Cnt}(P,Q)$ is the right inverse.
For $y\in X$, the fiber $\sigma _{s,y}$ of $\sigma _s$ is of the form $\sigma ^0_s\otimes \mathrm{id}_{B(l^2)}$.
For simplicity, in the rest of this subsection, we denote $\sigma ^0_s\otimes \mathrm{id}_{B(l^2)}$ by $\sigma _s$, if there are no danger of confusion.
Since $n_k$ is of the form $n'_k \otimes \mathrm{id}_{B(l^2)}$ for some $n'_k\in \mathrm{Aut}(P,Q)$, $n_k$ extends to an automorphism $\tilde{n}_k$ of $(M_y\supset N_y)\rtimes _\sigma \chi (P,Q)$. 
By the above proposition, for $\{ v'_n\}\subset \mathcal{U}(Q)$ with $\lim _{n\to \infty}\mathrm{Ad}v'_n =n'_k$, we have $\mathrm{Ad}(v'_n\otimes 1 )\to \tilde{n}_k$. 
Since $\tilde{n}_k$ commutes with the dual action $\hat{\sigma }$, the automorphism $\tilde{n}_k$ extends to $(M_y\supset N_y)\rtimes _{\sigma }\chi (P,Q)\rtimes _{\hat{\sigma }}\hat{\chi}(P,Q)$ by $\tilde{\tilde{n}}_k(\lambda _p)=\lambda _p$ for $p\in \hat{\chi} (P,Q)$. 
By using the above proposition again, we have $\mathrm{Ad}(v'_n\otimes 1) \to \tilde{\tilde{n}}_k$.
By Theorem \ref{masuda}, the map $s\mapsto u(n_k, \sigma _s)$ is $\sigma $-cocycle.
 Since $\sigma $ is a right inverse of the quotient, it is outer, as an action on $N_k$.
 Hence it is minimal because the group $\chi (P,Q)$ is a finite group.
 Hence any cocycle is a coboundary. 
Hence there exists a unitary $w$ of $N_y$ with $u(n_k, \sigma _s)=w^*\sigma (w)$. 
The automorphism $\mathrm{Ad}w \circ \tilde{\tilde{n}}_k$ is of the following form. 

\[ M_y\ni x \to \mathrm{Ad}(w)\circ n_k(x)\in M_y,\]
\[ \lambda _s\mapsto \lambda _s\]
for $s\in \chi (P,Q)$,
\[ \lambda _p \mapsto \lambda _p\]
for $p\in \hat{\chi} (P,Q)$.
Note that $\mathrm{Ad}w\circ n$ and $\sigma $ commute. 
Hence by Lemma 13.3 of Haagerup--St\o rmer \cite{HS} (See also Lemma 13.2 of Haagerup--St\o rmer), if we identify $(M_y\supset N_y) \rtimes _{\sigma }\chi (P,Q)\rtimes _{\hat{\sigma}}\hat{\chi}(P,Q)$ with $(M_y\supset N_y) \otimes B(l^2(\chi (P,Q)))$, there exists an automorphism $F$ of $(M_y\supset N_y) \otimes  B(l^2(\chi (P,Q)))$ such that 
\[ \mathrm{Ad}w \circ \tilde{\tilde{n}}_k =F\circ (\mathrm{Ad}w\circ n _k\otimes \mathrm{id}_{B(l^2(\chi (P,Q)))} )\circ F^{-1}.\]
We also have 
\begin{align*}
n_k\otimes \mathrm{id}_{B(l^2(\chi (P,Q)))}&\sim (\theta '\circ (n_k \otimes \mathrm{id}_{R_0})\circ (\theta ')^{-1})\otimes \mathrm{id}_{B(l^2(\chi (P,Q)))} \\
                                           &\sim (n_k \otimes \mathrm{id}_{R_0})\otimes \mathrm{id}_{B(l^2(\chi (P,Q)))} \\
                                           &\sim n_k \otimes \mathrm{id}_{R_0} \sim n_k,
\end{align*}
where $\alpha \sim \beta$ means that two actions $\alpha $ and $\beta $ are mutually conjugate.
Hence there exists an automorphism $\theta $ of $(M_y\supset N_y) \otimes B(l^2(\chi (P,Q)))$ satisfying
\[ \mathrm{Ad}w \circ \tilde{\tilde{n}}_k =\mathrm{Ad}(F(w\otimes 1_{B(l^2(\chi (P,Q)))})\theta (w^*)) \circ \theta \circ \mathrm{Ad} w\circ n_k \circ  \theta ^{-1}.\]
Set
\[\tilde{\tilde{m}}_{h,y}:=\theta \circ m_{h,y} \circ \theta ^{-1},\]
\begin{align*}
 \tilde{\tilde{\mathcal{A}_0}}(x):=\{ &(\rho ^x, w^x_h) \mid \rho ^x\in \mathrm{Aut} (M_x, N_x ),  \\ 
                                      &\{ w_h^x\} _{h\in H(\alpha )} \ \mathrm{is} \  \mathrm{an} \ \tilde{\tilde{m}}_x|_{H(\alpha)} \ \mathrm{cocycle}, \ \rho ^x \circ \tilde{\tilde{m}}_{h,x} \circ (\rho ^x)^{-1}=\mathrm{Ad}w^x_h \circ m_{h,x}\},
 \end{align*}
\[ \tilde{\tilde{\mathcal{B}_0}}(x):=\{ (\mathrm{Ad}u^x, u^x\tilde{\tilde{m}}_{h,x}(u^x)^*) \mid u^x \in \mathcal{U}(N_x) \}.\]

\begin{lemm}
\label{reduction}
Let $\theta$, $m_{h,y}$ and $\tilde{\tilde{m}}_{h,y}$ be as above.
Assume that automorphisms $\rho$, $\tilde{\tilde{\rho}}$ and a unitary $v $ of $N_y$ satisfy $\theta \circ \rho \circ \theta ^{-1}=\mathrm{Ad}v\circ \tilde{\tilde{\rho}}$. 
Assume that there exists an $m_{h,y}$-cocycle $\{ w_h\}$ with $[(\rho , w_h)]\in \mathcal{A}(x)$.  
Then $[(\rho ,  w_h)]\in \overline{\mathrm{Int}}B(y)$ if and only if $[(\tilde{\tilde{\rho}} , v^*\theta (w_h)\tilde{\tilde{m}}_{h,y}(v))]\in \overline{\mathrm{Int}}\tilde{\tilde{B}}(y)$.
\end{lemm} 
\begin{proof}
If $\{[(\mathrm{Ad}u^y_n, u^y_n\tilde{\tilde{m}}_{h,y}(u^y_n)^*)]\}$ converges to $[(\tilde{\tilde{\rho}} , v^*\theta (w_h)\tilde{\tilde{m}}_{h,y}(v))]$,
then 
\[ \{[(\mathrm{Ad}\theta ^{-1}(vu^y_n), \theta ^{-1}(vu^y_n)m_{h,y}(\theta ^{-1}(vu^y_n))^*)]\}\]
 converges to $[(\rho, w_h)]$. The converse is also true.
\end{proof}

\begin{lemm}
$[(\mathrm{Ad}w\circ \tilde{\tilde{ n}}, F(w\otimes 1_{B(l^2(\chi (P,Q)))})\theta (w_h^k)\tilde{\tilde{m}}_{h,y}(F(w\otimes 1_{B(l^2(\chi (P,Q)))}))^*]$ is approximately inner. 
\end{lemm}
\begin{proof}
Set 
\[ n'_k:=\mathrm{Ad}w \circ n_k,\]
\[ w'_h:= F(w\otimes 1_{B(l^2(\chi (P,Q)))})\theta (w^k_h)\tilde{\tilde{m}}_{h,y}(F(w\otimes 1_{B(l^2(\chi (P,Q)))}))^*.\]
Then the family $\{ w'_h\}$ is an $\tilde{\tilde{m}}_{n,y}$-cocycle and we have 
\[ \tilde{\tilde{n}}'_k \circ \tilde{\tilde{m}}_{h,y}\circ (\tilde{\tilde{n}}'_k)^{-1}=\mathrm{Ad}w'_h \circ \tilde{\tilde{m}}_{h,y}.\]
The following is a slight modification of the proof of Lemma 2.4 of Kawahigashi--Takesaki \cite{KT}. 
Since $n'_k$ preserves a trace and has the trivial Loi invariant, the automorphism $n'_k$ is approximately inner. 
Hence there exists a sequence of unitaries $\{ v_n\}$ of $N_y$ with $\lim _{n\to \infty } \mathrm{Ad}v_n=n'_k$.
 By (1) and (2) of Proposition \ref{3.2 of masuda}, we have 
\[ \lim \mathrm{Ad}v_n=\tilde{\tilde{n}}'_k.\]
 Note that $v_n \in N_y^{\hat{\hat{\sigma}}} $ and $\hat{\hat{\sigma}}$ commutes with $\tilde{\tilde{n}}'_k$. 
For $h\in N(\alpha)$, there exists a unitary $a^y_h\in N_y$ with
\[ \tilde{\tilde{m}}_{h,y}=\mathrm{Ad}(a^y_h)\circ \hat{\hat{\sigma}} _{\nu (h)}.\]
This is because $\tilde{\tilde{m}}_{h,y}(=\theta \circ m_{h,y} \circ \theta ^{-1})$ and $m_{h,y}$ are in the same class in $\chi (M_y,N_y)$ by Theorem \ref{masuda}.
Hence we have
\begin{align*}
\mathrm{Ad}(w'_h a^y_h)\hat{\hat{\sigma }}_{\nu (h)} &=\mathrm{Ad}w'_h\circ \tilde{\tilde{m}}_{h,y} \\
                                                            &=\tilde{\tilde{n}}'_k\circ \tilde{\tilde{m}}_{h,y} \circ (\tilde{\tilde{n}}'_k)^{-1} \\
                                                            &=\mathrm{Ad}(\tilde{\tilde{n}}'_k(a^y_h))\hat{\hat{\sigma }}_{\nu (h)}.
\end{align*}
Hence there exists $c^k_h\in \mathbf{T}$ with
\[ \tilde{\tilde{n}}'_k(a^y_h)=c^k_hw'_ha^y_h.\]
 Since $\mu (h,h')a_{hh'}^y=a_h^y\hat{\hat{\sigma }}_{\nu (h)}(a_{h'}^y)$ for some $\mu (h,h')\in \mathbf{T}$, we have

\begin{align*}
c^k_hc^k_{h'}\mu (h,h')w'_{hh'}a_{hh'}^y &=c^k_hc^k_{h'}w'_h\tilde{\tilde{m}}_{h,y}(w'_{h'}))a^y_h\hat{\hat{\sigma }}_{\nu (h)}(a^y_{h'}) \\
                                          &=c^k_hw'_h\tilde{\tilde{m}}_{h,y}(c^k_{h'}w'_{h'})a^y_h\hat{\hat{\sigma}}_{\nu(h)}(a^y_{h'}) \\
                                          &=\tilde{\tilde{n}}'_k(a^y_h)\hat{\hat{\sigma }}_{\nu (h)}(\tilde{\tilde{n}}'_k(a^y_{h'})) \\
                                          &=\tilde{\tilde{n}}'_k(a^y_h\hat{\hat{\sigma}}_{\nu (h)}(a^y_{h'})) \\
                                          &=\tilde{\tilde{n}}'_k(\mu (h,h')a_{hh'}^y) \\
                                          &=\mu (h,h')c^y_{hh'}w'_{hh'}a^y_{hh'}. 
\end{align*}
In the fourth equality, we used the fact that $\tilde{\tilde{n}}'_k$ commutes with $\hat{\hat{\sigma }}_{\nu(h)}$. 
Hence the map $N(\alpha )\ni h\mapsto c^k_h$ is a character. 
This $c^k$ is extended to a character $\overline{c}^k$ of $H(\alpha )$. 
Since $[(\tilde{\tilde{n}}'_k, w'_h)]=[(\tilde{\tilde{n}}'_k, \overline{c}^k_hw'_h))]$, we may assume that $\overline{c}^k=1$.
Now, we have
\[ \tilde{\tilde{n}}'_k(a^y_h)=w'_ha^y_h,\]
\[ \tilde{\tilde{n}}'_k\circ \tilde{\tilde{m}}_{h,y}\circ (\tilde{\tilde{n}}'_k)^{-1}=\mathrm{Ad}(w'_h)\circ \tilde{\tilde{m}}_{h,y}\]
and $\tilde{\tilde{n}}'_k=\lim _{n\to \infty}\mathrm{Ad}v_n$.
Since 
\[ \mathrm{Ad}(w'_h)^*v_n\tilde{\tilde{m}}_{h,y}(v_h^*)) \to \mathrm{Ad}(w'_h)^*\tilde{\tilde{n}}'_k\circ \tilde{\tilde{m}}_{h,y}\circ (\tilde{\tilde{n}}'_k)^{-1}=\tilde{\tilde{m}}_{h,y}\]
and $\tilde{\tilde{m}}_{h,y}$ is centrally trivial, the sequence 
\[ \{ (w'_h)^*v_n\tilde{\tilde{m}}_{h,y}(v_n^*)\} \]
is centralizing. We also have
\begin{align*}
 (w'_h)^*v_n\tilde{\tilde{m}}_{h,y}(v_n^*) &=(w'_h)^*v_na^y_hv_n^*(a^y_h)^* \\
                                                   &\to (w'_h)^*\tilde{\tilde{n}}'_k(a^y_h)(a^y_h)^* =1.
\end{align*}
Thus rest of the proof is the same as that of Lemma 2.4 of Kawahigashi--Takesaki \cite{KT}.
\end{proof}

Now, we return to the proof of Lemma \ref{point}.

\bigskip

\textit{Proof.}
By the previous lemma, $[(\mathrm{Ad}w\circ \tilde{\tilde{n}}_k, w\theta (w_h^k)\tilde{\tilde{m}}_{h,y}(w^*))]$ is approximately inner. 
Hence by Lemma \ref{reduction}, $[(\mathrm{Ad}w\circ n, ww_h^km_{h,y}(w^*))]$ is approximately inner. Note that we applied Lemma \ref{reduction} for $(\rho , v, w_h)\leftrightarrow (\mathrm{Ad}w\circ n, F(w\otimes 1), ww_h^k m_{h,y}(w^*))$. Hence by Lemma \ref{inner}, $[(n,w^k_h)]$ is approximately inner. \qed

\begin{rema}
\label{end}
In the above, we have only considered the case when $\mathrm{Aut}(P,Q)=\overline{\mathrm{Int}}(P,Q)$. We need to consider other cases in Remark \ref{other}.
 However, other cases are shown by the same argument as that of the proof of the case of single factors (Theorem 2.1.14 of Jones--Takesaki \cite{JT}). 
The following two observations are crucial.

(1) Actions of discrete abelian groups on inclusions of semifinite factors satisfying $\mathrm{Int}(P,Q)=\mathrm{Cnt}(P,Q)$ have been classified, up to cocycle conjugacy (Theorem 5.1 of Masuda \cite {M}). 
Actions of discrete abelian groups on inclusions of the form $P=Q\rtimes \mathbf{Z}/p\mathbf{Z}\supset Q$ are classified, up to cocycle conjugacy. 
This can be achieved by the same argument as that of Sections 3--6 of Masuda \cite{M2} (Sections 3, 4 and 6 are completely same. 
For Section 5, the important observation is that we have $\mathcal{Z}(Q) \cap C_\omega (P,Q)=\mathbf{C}$. 
Hence it is possible to show a lemma similar to Lemma 5.2 of Masuda \cite{M2} (See also Lemma 5.2 of Masuda--Tomatsu \cite{MT})). 
Hence it is possible to show a lemma corresponding to Lemma \ref{factors}.

(2) Any centrally trivial automorphism of $P\supset Q$ is implemented by a unitary $u$ of $P$ or $P_1$, where $P_1$ is the basic extension of $P\supset Q$. 
Hence it is possible to show the statement corresponding to Lemma \ref{point} by the same argument as that of the proof of Lemma 2.5.6 of Jones--Takesaki \cite{JT}.
\end{rema}  
 
 \section{Classification of actions of compact abelian groups}
 \label{compact}
 \subsection{Classification up to stable conjugacy}
In this subsection, we show a classification theorem of actions of compact abelian groups on any amenable subfactor with index less than $4$, up to stable conjugacy. 
\begin{defi}
\textup{(Jones--Takesaki \cite{JT})} Let $\alpha $ and $\beta $ be two actions of a locally compact group $G$ on an inclusion $M\supset N$. 
The action $\alpha $ is said to be stably conjugate to $\beta$ if $\alpha \otimes \mathrm{Ad}\rho $ is conjugate to $\beta \otimes \mathrm{Ad}\rho $, where, $\rho $ is the right translation of $L^2G$.
\end{defi}
\begin{theo}
\label{cptabl}
Let $M \supset N$ be an inclusion of AFD factors of type $\mathrm{II}_1 $ with index less than $4$ and $A$ be a compact abelian group. 
Let $\alpha $ and $\beta $ be two actions of $A$ on the inclusion $M\supset N$. 
Then they are mutually stably conjugate if and only if there exists an isomorphism $\theta $ from $(M\supset N) \rtimes _\alpha A$ to $(M \supset N) \rtimes _\beta A$ which satisfies the following conditions.

\textup{(1)} For any $g\in \hat{A}$, we have $\Phi (\theta \circ \hat{\alpha }_g\circ \theta ^{-1})=\Phi (\hat{\beta}_g)$.

\textup{(2)} We have $N(\hat{\alpha})=N(\hat{\beta} )$.

\textup{(3)} We have $\nu _{\theta \circ \hat{\alpha} \circ \theta ^{-1}}=\nu _{\hat{\beta}}$.

\textup{(4)} We have $[\theta (\lambda ^{\hat{\alpha }}, \nu ^{\hat{\alpha}})]=[(\lambda ^{\hat{\beta}},\mu ^{\hat{\beta}})]$.
\end{theo}
\begin{proof}
By Proposition \ref{probability},  the crossed product $(M\supset N) \rtimes _\alpha A$ has the same probability index as that of the original inclusion $M\supset N$.
By the central ergodicity of $\hat{\alpha }$ and Theorem \ref{ergodicity}, the crossed product inclusion is isomorphic to an inclusion of the form $A\otimes B(l^2) \otimes (P\supset Q)$, where $A$ is an abelian von Neumann algebra and $B(l^2)$ is a type I factor and $P\supset Q$ is an inclusion satisfying one of the following conditions.

(1) The inclusion $P\supset Q $ is a subfactor of type $\mathrm{II}_1$ such that we have $\mathrm{Aut}(P,Q)=\overline{\mathrm{Int}}(P,Q)$ and there exists a unique right inverse $\sigma^0:\chi (P,Q)\to \mathrm{Cnt}(P,Q)$, up to (cocycle) conjugacy.

(2) The inclusion $P\supset Q $ is a subfactor of type $\mathrm{II}_1$ satisfying $\mathrm{Cnt}(P,Q)=\mathrm{Int}(P,Q)$.

(3) The inclusion $P\supset Q $ is a subfactor of type $\mathrm{II}_1$ with its principal graph is of type $D_4$.

(4) $Q$ is a finite factor and $P=Q\otimes \mathbf{C}^p$, where $p=2$ or $3$.

(5) $Q=R\otimes \mathbf{C}^p$, where $R$ is a finite factor and $p=2$ or $3$, $P=R\otimes M_p(\mathbf{C})$.

Hence it is possible to apply Theorem \ref{classify}  and Remark \ref{end} of the previous section. 
Hence by the same argument as that of the proof of Theorem 3.1.3 of Jones--Takesaki \cite{JT} or Theorem 3.1 of Kawahigashi--Takesaki \cite{KT}, we are done.
\end{proof}

\subsection{Classification up to conjugacy}
In this subsection, we explain a classification theorem up to conjugacy. 
The following argument is completely the same as that of Subsection 3.2 of Jones--Takesaki \cite{JT}.
Let $A$ be a compact abelian group and $M \supset N$ be a subfactor with a normal faithful conditional expectation $E:M\to N$. 
Let $\mathcal{S}$ be the set of all stable conjugacy classes of actions of $A$ on $M \supset N$ commuting with the expectation $E$.
For an action $\alpha $ of $A$ on $M \supset N$, we denote the stable conjugacy class of $\alpha$ by $\partial (\alpha )$.
For $s\in \mathcal{S}$, we choose an action $m $ of $A$ satisfying $\partial (m ) =s$.
Set $p:=\int _A 1\otimes \rho _t \ d\mu (t) \in N \otimes B(L^2(A))$, where $\rho _t$ is the regular representation of $A$ and $\mu $ is the Haar measure of $A$.
Let $\mathcal{H}_{s,M\supset N}$ be the set of all projections of $N \rtimes _m A$($=N\otimes B(L^2A)^{m \otimes \mathrm{Ad}\rho}$).
We denote the set of all automorphism of $(M \supset N) \otimes B(L^2A)$ which commutes with $m \otimes \mathrm{Ad}\rho$ by  $\mathrm{Aut}_{m\otimes \mathrm{Ad}\rho }(M,N)$.
Then $\mathrm{Aut}_{m\otimes \mathrm{Ad}\rho }(M,N)$ acts on $\mathcal{H}_{s,M\supset N}$. 
Let $\mathcal{P}_{s, M \supset N}$ be the set of all orbits of the action of $\mathrm{Aut}_{m\otimes \mathrm{Ad}\rho }(M,N)$ on $\mathcal{H}_{s,M\supset N}$.
For any action $\alpha $ with $\partial (\alpha )=s$, it is possible to choose an automorphism $\theta \in \mathrm{Aut}((M\supset N) \otimes B(L^2A))$ with $\theta \circ (\alpha \otimes \mathrm{Ad}\rho )\circ \theta ^{-1} =m\otimes \mathrm{Ad}\rho$.
Then as shown in Subsection 3.2 of Jones--Takesaki \cite{JT}, $\theta (p) \in \mathcal{H}_{s, M \supset N}$.
The orbit of $\theta (p) $ in $\mathcal{P}_{s, M \supset N}$ does not depend on the choice of $\theta$.
Let $\iota (\alpha )$ be the orbit of $\theta (p)$ in $P_{s, M \supset N}$, which is said to be the inner invariant of $\alpha$.
We have the following proposition.
\begin{prop}
Let $\alpha $ and $\beta$ be tow actions of $A$ on $M \supset N$. 
Then $\alpha $ and $\beta $ are mutually conjugate if and only if $\partial(\alpha)=\partial(\beta)$ and $\iota (\alpha )=\iota (\beta)$.
\end{prop} 
\begin{proof}
This is shown by completely the same argument as that of the proof of Proposition 3.2.2 of Jones--Takesaki \cite{JT}.
\end{proof}

\section{Actions of $\mathbf{T}$}
In many cases, invariants of Theorem \ref{cptabl} are difficult to compute.
Although Theorem \ref{cptabl} provides a list of obstructions for stable conjugacy, it may not be better to think that it is useful for actual computation.
 However, for an action $\alpha $ of the one dimensional torus $\mathbf{T}$, the invariants are turned out to be easy because $N(\hat{\alpha} )$ is trivial.

\begin{theo}
\label{less 4}
Let $M \supset N$ be an inclusion of AFD factors of type $\mathrm{II}_1$ with index less than $4$. 
Then the actions $\alpha$ of $\mathbf{T}$ on the inclusion $M\supset N$ are classified by the conjugacy crass of the Loi invariant $\Phi(\hat{\alpha})$ of the dual actions $\hat{\alpha}$, up to stably conjugacy. 
The actions with Connes spectra non-zero are classified by their crossed products, up to cocycle conjugacy. 
The list of the crossed products of actions with Connes spectra $L\mathbf{Z}$ \textup{(}$L\not =0$\textup{)}, as actions on $N$ is the following.

\begin{table}[htb]
  \begin{tabular}{|l|c|} \hline
    $M\supset N$ \textup{(}type $\mathrm{II}_1$\textup{)} & Crossed Product \textup{(}type $\mathrm{II}_\infty$\textup{)} \\ \hline \hline
    $A_3$ & $A_3 \otimes \mathbf{C}^L$ \\
          & $(\mathbf{C}^2\supset \mathbf{C}) \otimes R_{0,1} \otimes \mathbf{C}^L$ \\
          & $(M_2(\mathbf{C}) \supset \mathbf{C}^2) \otimes R_{0,1} \otimes \mathbf{C}^{L/2}$ \textup{(} \textup{if} $L\in 2\mathbf{Z}$\textup{)} \\ \hline
   $A_m$ \textup{(}$m=4n-3$, $n\geq 2$\textup{)} & $D_{2n}\otimes \mathbf{C}^L$  \\ 
          & $A_m  \otimes \mathbf{C}^L$    \\ \hline
   $A_m$ \textup{(}$m\not =4n-3,3$\textup{)}& $A_m \otimes \mathbf{C}^L$ \\ \hline
   $D_4$  & $D_4\otimes  \mathbf{C}^L$ \\
          & $(\mathbf{C}^3\supset \mathbf{C}) \otimes R_{0,1} \otimes \mathbf{C}^L$ \\
          & $(M_3(\mathbf{C})\supset \mathbf{C}^3) \otimes R_{0,1} \otimes \mathbf{C}^{L/3}$ \textup{(} $L\in 3\mathbf{Z}$\textup{)} \\ \hline
   $D_{2n}$ \textup{(}$n\geq 3$\textup{)} & $D_{2n} \otimes \mathbf{C}^L$  \\ \hline 
   $E_6$  & $E_6  \otimes \mathbf{C}^L$ \textup{(itself, not the anti-isomorphic one)} \\ \hline
   $E_8$  & $E_8 \otimes  \mathbf{C}^L$ \textup{(itself, not the anti-isomorphic one)} \\ \hline

  \end{tabular}
\end{table}
Here, $A_m$, $D_{2n}$, $E_6$ and $E_8$ denote the inclusions of factors with their principal graphs of type $A_m$, $D_{2n}$, $E_6$ and $E_8$, respectively.
\end{theo}
In the following, we will show Theorem \ref{less 4}.
In order to achieve this, the following proposition, which is similar to Proposition 7 of Choda--Kosaki \cite{CK}, is important.

\begin{prop}
\textup{(See Proposition 7 of Choda--Kosaki \cite{CK}).}\label{strong outer} Let $M\supset N$ be a finite depth inclusion of factors of type $\mathrm{II}_1$ and $\alpha $ be an action of $\mathbf{T}$ on $M\supset N$.
Set $\tilde{M} \supset \tilde{N} :=(M\supset N) \rtimes _\alpha \mathbf{T}$.
Assume that we have $\mathcal{Z}(M)=\mathcal{Z}(N)$.
Let $\tilde{E}:\tilde{M} \to \tilde{N}$ be the expectation with respect to a trace.  
Fix the Jones tower 
\[ \tilde{N} \subset \tilde{M}\subset \tilde{M_1}\subset \tilde{M_2} \subset \cdots\]
of the inclusion $\tilde{M}\supset \tilde{N}$. 
Then we have
\[ (\tilde{M}_k \cap \tilde{N}')^{\hat{\alpha}}=(\tilde{M}_k \rtimes _{\hat{\alpha}}\mathbf{Z})\cap (\tilde{N}\rtimes _{\hat{\alpha}}\mathbf{Z})'.\]
\end{prop}
\begin{proof}
The inclusion ``$\subset $ '' is trivial. 
We show the opposite direction. 
Choose $x\in (\tilde{M}_k \rtimes _{\hat{\alpha}}\mathbf{Z})\cap (\tilde{N} \rtimes _{\hat{\alpha}}\mathbf{Z})'$. 
Let $x=\sum _{n\in \mathbf{Z}}x_n\lambda ^{\hat{\alpha}}_n$ be the Fourier expansion. 
Since $x$ commutes with any $y\in \tilde{N}$, we have $yx_n =x_n \hat{\alpha}_n(y)$ for any $y\in \tilde{N}$, $n\in \mathbf{Z}$. 

\textbf{Case 1.} The case when $\hat{\alpha}|_{\mathcal{Z}(\tilde{M})}$ is free.

Assume that $x_n\not =0$ for some non-zero $n$. 
Since $\hat{\alpha }|_{\mathcal{Z}(\tilde{M})}$ is free, there would exist a non-zero central projection $z$ such that $\hat{\alpha}_n(z)z=0$, and $z\leq (\mathrm{central} \ \mathrm{support} \ \mathrm{of} \ x_n)$. 
Then we have
\begin{align*}
|zx_n|^2&=x_n^* z^2 x_n \\
        &=x_n^* z x_n \hat{\alpha}_n (z) \\
        &=|x_n|^2z\hat{\alpha }_n(z) =0.
\end{align*}
On the other hand, since $z\leq (\mathrm{central} \ \mathrm{support} \ \mathrm{of} \ x_n)$,  we would have $zx_n\not = 0$, which would be a contradiction. Hence we have $x_n =0$ for all $n\not =0$.

\textbf{Case 2.} The case when $\hat{\alpha}|_{\mathcal{Z}(\tilde{M})}$ has a non-zero period $L$.

In this case we have $\mathcal{Z}(\tilde{M})=\oplus _{l=1}^L\mathbf{C}p_l$ for the minimal central projections $\{ p_l\}$.
 We show that for each $l=1, \cdots, L$, the family $\{ \hat{\alpha }_{Ln}\}_{n\in \mathbf{Z}}$ defines a strongly outer action of $\tilde{M}_{p_l}\supset \tilde{N}_{p_l}$. 
This is shown by the following way.
Assume that $\{ \hat{\alpha }_{Ln}\}_{n\in \mathbf{Z}}$ were not strongly outer. 
Then by Corollary 4 of Choda--Kosaki \cite{CK}, there would exist a positive integer $n>0$ and a unitary $u\in \mathcal{U}(\tilde{N}_{p_1})$ with $\hat{\alpha }_{Ln}|_{\tilde{M}_{p_1}}=\mathrm{Ad}(u)$. 
Since $\hat{\alpha }_{Ln} =\mathrm{Ad}(\hat{\alpha}_L(u))$, there would exist a complex number $\gamma $ with $|\gamma |=1$, $\gamma ^n =1$ and $\hat{\alpha }(u)=\gamma u$ (Connes obstruction). 
Then we would have $\hat{\alpha }_L(u^n) =u^n$ and $\mathrm{Ad}((u^{n})^*)\circ \hat{\alpha }_{Ln^2}|_{\tilde{M}_{p_1}}=\mathrm{id}_{\tilde{M}_{p_1}}$. 
Set
\[ v:=\sum _{l=1}^L \hat{\alpha}_{l-1}((u^n)^*).\]
Then $v$ is a unitary of $\tilde{N}$.
 Let $\{ \lambda ^{\hat{\alpha}}_{n} \}_{n\in \mathbf{Z}}$ be the unitaries of $\tilde{M}\rtimes _{\hat{\alpha }}\mathbf{Z}$ satisfying $\lambda ^{\hat{\alpha}}_n x\lambda ^{\hat{\alpha}}_{-n} =\hat{\alpha}_n(x)$ for $x \in \tilde{M}$. 
Then we have $0\not =v\lambda ^{\hat{\alpha}}_{Ln^2}\in \mathcal{Z}(\tilde{M} \rtimes _{\hat{\alpha}}\mathbf{Z})$, which would contradict to the fact that $\tilde{M} \rtimes _{\hat{\alpha}}\mathbf{Z}\cong M\otimes B(L^2\mathbf{T})$ is a factor.
Hence for each $l=1, \cdots, L$, the family $\{ \hat{\alpha }_{Ln}\}_{n\in \mathbf{Z}}$ defines a strongly outer action of $\tilde{M}_{p_l}\supset \tilde{N}_{p_l}$.  
Hence we have $x_{Ln}=0$ for any non-zero $n\not =0$. 
On the other hand, by the same argument as that in Case 1, we have $x_j=0$ for $j\not \in L\mathbf{Z}$.
\end{proof}

\begin{lemm}
Let $M \supset N$ be an inclusion of AFD factors of type $\mathrm{II}_1$ with index less than $4$. 
Then the actions $\alpha$ of $\mathbf{T}$ on the inclusion $M\supset N$ are classified by the conjugacy crass of the Loi invariant $\Phi(\hat{\alpha})$ of the dual actions $\hat{\alpha}$, up to stably conjugacy. 
\end{lemm} 
\begin{proof}
By the proof of Proposition \ref{strong outer}, for any action $\alpha $ of $\mathbf{T}$ on $M\supset N$, we have $N(\hat{\alpha})=\{0\}$.
Hence by Theorem \ref{cptabl}, we get the conclusion.
\end{proof}

By the same argument as that of Proposition 6.2 of Masuda--Tomatsu \cite{MT}, it is possible to show that stable conjugacy implies cocycle conjugacy under the assumption that Connes spectra are non-zero.
Hence what we need to do in the following is to show that the  list of the crossed products in Theorem \ref{less 4} is correct.

 \begin{prop}
 \label{graph change}
Let $M\supset N $ be an inclusion of AFD factors of type $\mathrm{II}_1$ with index less than $4$ and $\alpha $ be an action of $\mathbf{T}$ on $M\supset N$.
Assume that the Connes spectrum of $\alpha $ is not zero.
Set $\tilde{M} \supset \tilde{N} :=(M\supset N) \rtimes _\alpha \mathbf{T}$.
Assume that $\mathcal{Z}(\tilde{M})=\mathcal{Z}(\tilde{N})$.
Then the following two conditions are equivalent.

\textup{(a)} The inclusion $M\supset N$ is isomorphic to any direct summand $P\supset Q$ of the inclusion $\tilde{M}\supset \tilde{N}$.

\textup{(b)} The action $\alpha $ is cocycle conjugate to an action of the form $\mathrm{id}_{M\supset N} \otimes \sigma $ for an action $\sigma $ on the AFD factor $R_0$ of type $\mathrm{II}_1$. 
\end{prop}

\begin{proof}
We denote $M_k \rtimes _\alpha \mathbf{T}$ by $\tilde{M}_k$.
Then by Proposition 3.13 of Kosaki \cite{Ko}, the sequence $\{ \tilde{M}_k\}$ is the Jones tower of the expectation $\tilde{E}:\tilde{M}\to \tilde{N}$ with respect to a trace.
By Popa's classification theorem (Theorem 5.1.1 of Popa \cite{P}), the implication (b) $\Rightarrow $ (a) is trivial. 
We show the opposite direction. 
It is enough to show that when the standard invariant of $M\supset N$ is isomorphic to that of $P\supset Q$, then condition (b) is satisfied. 
Let $\{ p_l\}_{l=1}^L$ be the minimal projections of $\mathcal{Z}(\tilde{M})=\mathcal{Z}(\tilde{N})$.
There exists a natural embedding 
\[ (\tilde{M_k}\cap \tilde{N}')^{\hat{\alpha}}p_l\subset (\tilde{M}_k \cap \tilde{N}')p_l =P_k \cap Q'.\]
They are finite dimensional. 
By condition (a) and Proposition \ref{strong outer}, we have $(\tilde{M_k}\cap \tilde{N}')^{\hat{\alpha}} p_l\cong P_k \cap Q'$. 
Hence we have $(\tilde{M}_k\cap \tilde{N}')^{\hat{\alpha}}p_l=(\tilde{M}_k \cap \tilde{N}')p_l$. 
This means that the Loi invariant of $\hat{\alpha}_L$ is trivial. 
Hence the Loi invariant of $\hat{\alpha}$ is conjugate to that of the following action.
\[ \mathrm{id}_{P\supset Q} \otimes \hat{\beta} \in \mathrm{Aut}((P\supset Q)\otimes (R_0\rtimes _\beta \mathbf{T})).\]
Here, $\beta$ is the (unique) action of $\mathbf{T}$ on $R_0$ with Connes spectrum $L\mathbf{Z}$. 
Hence by the above lemma, the action $\alpha$ is cocycle conjugate to $\beta$.
\end{proof}

\begin{lemm}
\label{nonfactor}
Let $M\supset N$ be an inclusion of AFD factors of type $\mathrm{II}_1$ and $\alpha $ be an action of $\mathbf{T}$ on the inclusion $M\supset N$. 
Assume that the inclusion $\tilde{M}\supset \tilde{N}:=(M\supset N) \rtimes _\alpha \mathbf{T}$ is of the form $\mathbf{C}^L\otimes (\mathbf{C}^p\supset \mathbf{C})\otimes B(l^2)\otimes R$ or $\mathbf{C}^L\otimes (M_p(\mathbf{C})\supset \mathbf{C}^p) \otimes B(l^2) \otimes R$, where $R$ is a finite factor, $L\in \mathbf{N}$ and $p$ is a prime number. 
Then there exists an action $\beta $ of $\mathbf{Z}/p\mathbf{Z}$ on $N$ satisfying $M=N\rtimes _\beta \mathbf{Z}/p\mathbf{Z}$.
\end{lemm}
\begin{proof}
We consider the case when $\tilde{M}\supset \tilde{N}:=(M\supset N) \rtimes _\alpha \mathbf{T}$ is of the form $\mathbf{C}^L\otimes (\mathbf{C}^p\supset \mathbf{C})\otimes B(l^2)\otimes R$. 
The other case is shown by a similar way.
Let $\lambda ^{\hat{\alpha }}_n$, $n\in \mathbf{Z}$ be an implementing unitary of $\hat{\alpha }$ and $u$ be a unitary of $\mathbf{C}^L \otimes \mathbf{C}^p \otimes \mathbf{C}\otimes \mathbf{C}$ be a unitary satisfying $\hat{\alpha }(u)=e^{2\pi i/p}u$. 
Note that this $u$ exists because $\hat{\alpha }$ acts ergodically on $\mathbf{C}^L\otimes \mathbf{C}^p \otimes \mathbf{C}\otimes \mathbf{C}$. 
Then we have $\lambda ^{\hat{\alpha }}_1u=e^{2\pi i/p}u\lambda ^{\hat{\alpha }}_1$. 
Hence $\mathrm{Ad}u$ preserves $\tilde {N} \rtimes _{\hat{\alpha}}\mathbf{Z}$ globally. 
We also have $u^p=1$ and that $\mathrm{Ad}u$ preserves the trace of $\tilde{N}\rtimes _{\hat{\alpha}}\mathbf{Z}$. Hence we have
\[ (\tilde{N}\rtimes _{\hat{\alpha }}\mathbf{Z})\rtimes _{\mathrm{Ad}u}(\mathbf{Z}/p\mathbf{Z})=\tilde{M}\rtimes _{\hat{\alpha}}\mathbf{Z}.\]
Since $\tilde{M}\rtimes _{\hat{\alpha}}\mathbf{Z}$ is a factor, $p$ is a prime number, by the uniqueness of the outer action of $\mathbf{Z}/p\mathbf{Z}$, the action $\mathrm{Ad}u$ is cocycle conjugate to an action of the form $\beta \otimes \mathrm{id}_{B(L^2\mathbf{T})}$, where $\beta $ is an outer action of $\mathbf{Z}/p\mathbf{Z}$ on $N$, by an identification through Takesaki's duality.
 Hence we get the conclusion. 
\end{proof}

\begin{lemm}
\label{D}
Let $M\supset N$ be an inclusion of AFD factors of type $\mathrm{II}_1$ with index less than $4$ and $\alpha $ be an action of $\mathbf{T}$ on the inclusion $M\supset N$. 
Assume that the inclusion $\tilde{M}\supset \tilde{N}$ is isomorphic to $\mathbf{C}^L\otimes D_{2n}$. Then the principal graph of the inclusion $M\supset N$ is of type $A_{4n-3}$ or of type $D_{2n}$.
\end{lemm}
\begin{proof}
The following is argument is essentially the same as that of the proof of Proposition 6.2 of Loi \cite{L}.
We need to exclude the possibility that the principal graph of $M\supset N$ is of type $E_6$ or $E_8$.
Let $k$ be the depth of $\tilde{M} \supset \tilde{N}$. 
Then the central support of the Jones projection $e_k$ in $\tilde{M}_k\cap \tilde{N}'$ is $1$. 
By Proposition \ref{strong outer}, the central support of $e_k$ in $M_k\cap N'$ is also $1$. 
Hence the depth of  $M\supset N$ is not greater than $k$.
On the other hand, by Proposition \ref{probability}, the index of $M\supset N$ and that of any direct summand of $\tilde{M}\supset \tilde{N}$ are same.
Recall that when the index is less than $4$, the index is given by the norm $4\cos ^2(\pi /h)$ of the principal graph, where $h$ is the Coxeter number of the principal graph ($h=m+1$ for $A_m$, $h=4m-2$ for $D_{2m}$, $12$ and $30$ for $E_6$ and $E_8$, respectively). 
Comparing indices and depth of $M\supset N$ and direct summands of $\tilde{M}\supset \tilde{N}$, if the principal graph of $M\supset N$ were of type $E_8$, then we would have $15=2m-1\geq 6$, which would be a contradiction.
The case $E_6$ is excluded in a similar way. 
\end{proof}

\textit{Proof of Theorem \ref{less 4}}.
Let $\alpha $ be an action of $\mathbf{T}$ on $M\supset N$.
Set 
\[ \tilde{M} \supset \tilde{N} := M\rtimes _\alpha \mathbf{T}\supset N \rtimes _\alpha \mathbf{T}.\]
Assume that Connes spectrum of $\alpha $ (as an action on $N$) is $L\mathbf{Z}$.

\textbf{Case 1}. Assume that the principal graph of a direct summand of $\tilde{M}\supset \tilde{N}$ is of type $D_{2n}$.
By Remarks 6.5 of Loi \cite{L}, the number of candidate for cocycle conjugacy classes of $\hat{\alpha }$ is at most $2$. 
By Lemma \ref{D}, the principal graph of $M\supset N$ should be of type $D_{2n}$ or of type $A_{4n-3}$.

\textbf{Case 2}. Assume that the center $\mathcal{Z}(N)$ of $N$ does not coincide with $\mathcal{Z}(M)$. 
Then by Lemma \ref{nonfactor}, the principal graph of $M\supset N$ is of type $D_4$ or $A_3$.

\textbf{Case 3}. Assume that the principal graph of a direct summand of $\tilde{M}\supset \tilde{N}$ is of type $A_m$, $E_6$ or $E_8$.
By lines 1--3 of p.288 of Loi \cite{L}, and Proposition \ref{graph change}, $(M\supset N)\otimes B(l^2)$ is isomorphic to a direct summand of $\tilde{M}\supset \tilde{N}$. 
\qed

\end{document}